\documentclass[a4paper,12pt]{article}

\usepackage[utf8]{inputenc}
\usepackage[T1]{fontenc}
\usepackage[left=2.5cm,right=2.5cm,top=1.5cm,bottom=2cm]{geometry}
\usepackage{amsmath,amssymb,makeidx}
\usepackage[english]{babel}

\usepackage{hyperref}
\usepackage[language=english,backend=biber]{biblatex}
\usepackage{fourier}
\usepackage{graphics}
\usepackage{graphicx}
\usepackage{amsthm}
\usepackage{faktor}
\usepackage{thmtools}
\usepackage{enumitem}
\newlist{itemizeth}{itemize}{2}
\setlist[itemizeth]{label=\textbullet,noitemsep,topsep=0 mm}
\newlist{enumerateth}{enumerate}{2}
\setlist[enumerateth]{label*=\textbf{\alph*)},itemsep=0.2mm}
\setlist[itemize]{fullwidth,topsep=0 mm, leftmargin=0pt}
\setlist[enumerate]{label*=\textbf{\arabic*)},itemsep=0.2mm, leftmargin=0pt}
\declaretheoremstyle[
title=Proof,
numbered=no,
headfont=\normalfont\bfseries,
notefont=\bfseries, notebraces={}{},
postheadspace=17pt,
bodyfont=\normalfont,
headindent=0pt,
qed=\qedsymbol,
spacebelow=3,
]{demostyle}

\declaretheoremstyle[
headfont=\normalfont\bfseries,
notefont=\mdseries, notebraces={(}{)},
bodyfont=\normalfont,
thmbox=S
]{standard}

\declaretheoremstyle[
title=Theorem,
headfont=\normalfont\scshape,
notefont=\mdseries, notebraces={(}{)},
thmbox=S,
]{theo}

\declaretheoremstyle[
title=Theorem,
headfont=\scshape\bfseries,
notefont=\mdseries, notebraces={(}{)},
thmbox=L,
]{import}

\declaretheoremstyle[, 
numbered=no
]{rem}

\declaretheorem[title=Definition,style=standard]{defi}
\declaretheorem[title=Proposition,style=standard]{prop}
\declaretheorem[title=Lemma,style=standard]{lem}

\declaretheorem[style=theo]{Th}

\declaretheorem[title=Remark,style=rem]{rqu}

\declaretheorem[style=demostyle]{dem}

\makeatletter
\newcommand*{\house}[1]{%
   \mathord{%
     \mathpalette\@house{#1}%
   }%
}
\newcommand*{\@house}[2]{%
   \dimen@=\fontdimen8 %
       \ifx#1\scriptscriptstyle\scriptscriptfont
       \else\ifx#1\scriptstyle\scriptfont
       \else\textfont\fi\fi
       3 %
   \sbox0{%
     $#1%
       \vrule width\dimen@\relax
       \overline{%
         \kern2\dimen@
         \begingroup 
           #2%
         \endgroup
         \kern2\dimen@
       }%
       \vrule width\dimen@\relax
       \mathsurround=1.5\dimen@ 
     $%
   }%
   \ht0=\dimexpr\ht0-\dimen@\relax
   \dp0=\dimexpr\dp0+2\dimen@\relax
   \vbox{%
     \kern\dimen@ 
     \copy0 %
   }%
}

\newcommand{\R}{\mathbb{R}}
\newcommand{\N}{\mathbb{N}}

\newcommand{\Z}{\mathbb{Z}}
\newcommand{\Q}{\mathbb{Q}}
\newcommand{\C}{\mathbb{C}}
\newcommand{\K}{\mathbb{K}}
\newcommand{\Oal}{\mathcal{O}}
\newcommand{\Qbar}{\overline{\mathbb{Q}}}

\newcommand{\Spec}{\mathrm{Spec} \,}

\newcommand{\GL}{\mathrm{GL}}

\newcommand{\ord}{\mathrm{ord} \,}

\newcommand{\gf}{\dfrac}
\newcommand{\Vect}{\mathrm{Span}}
\newcommand{\Img}{\mathrm{Im} \,}
\newcommand{\den}{d}
\newcommand{\ddz}{\gf{\mathrm{d}}{\mathrm{d}z}}
\newcommand{\Sym}{\mathrm{Sym}}

\newcommand{\fonction}[5]{\begin{array}[t]{lrcl} 
#1: & #2 & \longrightarrow & #3 \\
    & #4 & \longmapsto & #5 \end{array}}

\title{Quantitative problems on the size of $G$-operators}
\author{Gabriel Lepetit}
\date{\today}
\makeindex
\begin{document}

\maketitle

\begin{abstract}
$G$-operators, a class of differential operators containing the differential operators of minimal order annihilating Siegel's $G$-functions, satisfy a condition of moderate growth called \emph{Galochkin condition}, encoded by a $p$-adic quantity, the size. Previous works of Chudnovsky, André and Dwork have provided inequalities between the size of a $G$ -operator and certain computable constants depending among others on its solutions. First, we recall André's idea to attach a notion of size to differential modules and detail his results on the behavior of the size relatively to the standard algebraic operations on the modules. This is the corner stone to prove a quantitative version of André's generalization of Chudnovsky's Theorem: for $f(z)=\sum_{\alpha, k,\ell} c_{\alpha, k,\ell} z^{\alpha} \log(z)^k f_{\alpha, k,\ell}(z)$, where $f_{\alpha, k,\ell}(z)$ are $G$-functions, we can determine an upper bound on the size of the minimal operator $L$ over $\Qbar(z)$ of $f(z)$ in terms of quantities depending on the $f_{\alpha, k,\ell}(z)$, the rationals $\alpha$ and the integers $k$. We give two applications of this result: we estimate the size of a product of two $G$-operators in function of the size of each operator; we also compute a constant appearing in a Diophantine problem encountered by the author.
\end{abstract}

\textbf{Remarks and notations}
\begin{itemize}
    \item In this paper, given some function $f$, we will call "minimal operator of $f$ over $\Qbar(z)$", or "minimal operator" where there is no possible ambiguity, any nonzero operator $L \in \Qbar(z)\left[\mathrm{d}/\mathrm{d}z \right]$ such that $L(f(z))=0$, and whose order is minimal for $f$.

\item For $u_1, \dots, u_n \in \Qbar$, we denote by $\den(u_1, \dots, u_n)$ the \emph{denominator} of $u_1, \dots, u_n$, that is to say the smallest $d \in \N^*$ such that $d u_1, \dots, d u_n$ are algebraic integers.

\item If $\K$ is a subfield of $\Qbar$, we denote by $\Oal_{\K}$ the set of algebraic integers of $\K$.
\end{itemize}

\section{Introduction}

A \emph{$G$-function} is a power series $f(z)=\sum\limits_{n=0}^{\infty} a_n z^n \in \Qbar\llbracket z\rrbracket$ satisfying the following hypotheses:
\begin{enumerateth}
\item $f$ is solution of a nonzero linear differential equation with coefficients in $\Qbar(z)$;
\item There exists  $C_1 >0$ such that $\forall n \in \N, \; \house{a_n} \leqslant C_1^{n+1}$,  where, for $\alpha \in \Qbar$, $\house{\alpha}:=\max\limits_{\sigma \in \mathrm{Gal}(\Qbar/\Q)} |\sigma(\alpha)|$ is the \emph{house} of $\alpha$.
\item There exists $C_2 >0$ such that $\forall n \in \N, \;  \den(a_0, \dots,a_n) \leqslant C_2^{n+1}$.
\end{enumerateth}

This family of special functions has been studied together with the family of $E$-functions, which are the functions $f(z)=\sum\limits_{n=0}^{\infty} (a_n/n!) z^n$ satisfying \textbf{a)} and such that the $a_n$ satisfy the conditions \textbf{b)} and \textbf{c)}. Siegel defined both classes in \cite{Siegelarticle}. The most basic example of $G$-function, which gives it its name, is the geometric series $f(z)=\sum\limits_{n=0}^{\infty} z^n=1/(1-z)$. Other examples includes the polylogarithms functions $\mathrm{Li}_s(z):=\sum\limits_{n=1}^{\infty} z^n/n^s$, or some hypergeometric series with rational parameters.

The theory of $G$-functions was largely developed from the 1970s with the works of Galochkin \cite{Galochkin74}, Chudnovsky \cite{Chudnovsky} and then André's book \cite{Andre}, which is the first extensive and systematic review of this theory, later enhanced by Dwork \cite{Dwork}.

The corner stone of this theory is the study of the nonzero minimal operators of $G$-functions over $\Qbar(z)$. It turns out that they have specific properties, including the fact that they are Fuchsian differential operators with rational exponents. They belong to the class of  $G$-operators, which are believed to "come from geometry" (see \cite[chapter II]{Andre}).

Chudnovsky's Theorem states that the nonzero minimal operator $L$ of a nonzero $G$-function $f$ satisfies an arithmetic condition of geometric growth on some denominators, called \emph{Galochkin condition}. This condition can be equivalently expressed by the fact that a quantity $\sigma(L)$, which is defined in $p$-adic terms (see Definition \ref{def:taillematrice} in Section \ref{sec:sizediffmodule} below), is finite. This quantity is called the \emph{size} of $L$. More generally, it is possible to associate a size $\sigma(G)$ with a differential system $y'=Gy$, $G \in \mathcal{M}_n(\Qbar(z))$.

Similarly, we can associate with a $G$-function its size, encoding the conditions \textbf{b)} and \textbf{c)} above.

Chudnovsky's Theorem was proved in such a way that there is an explicit relation between the size of $G$ and the size of $f$ (see \cite[chapter VII]{Dwork}). This is particularly useful for Diophantine applications, as shown in \cite{Lepetit2}.
\bigskip

In \cite{AndregevreyI} and \cite{AndregevreyII}, André studied the properties of the \emph{Nilsson-Gevrey series of arithmetic type} of order $0$, \emph{i.e.} functions that can be expressed as a finite sum \begin{equation}\label{eq:defNGseries}f(z)=\sum_{(\alpha, k, \ell) \in S}^{} c_{\alpha, k,\ell} z^{\alpha} \log(z)^k f_{\alpha, k,\ell}(z)\end{equation} where $S \subset \Q \times \N \times \N$, $c_{\alpha, k,\ell} \in \C^*$, and the $f_{\alpha, k,\ell}(z)$ satisfy the conditions \textbf{b)} and \textbf{c)} above. We say moreover that $f(z)$ is \emph{of holonomic type} if the $f_{k,\ell}(z)$ are holonomic, \emph{i.e.} satisfy the condition~\textbf{a)} \footnote{
In the remark p. 717 of \cite{AndregevreyI}, André seems to indicate that all holonomic Nilsson-Gevrey series of arithmetic type are of holonomic type but we don't see how to prove it with his argument. Indeed, it is stated that if $f(z)$ is of the form  \eqref{eq:defNGseries} where the $\alpha$ are in distinct classes modulo $\Z$ and the $f_{\alpha,k,\ell}(z)$ satisfy the conditions \textbf{b)} and \textbf{c)}, then the $f_{\alpha,k,\ell}(z)$ are all holonomic, which is false as the example of $0=y_0(z)-y_0(z)$, with $y_0(z)$ non holonomic, shows it.}. 

André proved in particular the following analogue of Chudnovsky's Theorem for the Nilsson-Gevrey series of arithmetic and holonomic type of order $0$: 
\begin{Th}[\cite{AndregevreyI}, p. 720]\label{th:chudnovskynilssongevreyandreI}
Let $S \subset \Q \times \N \times \N$ be a finite set, $(c_{\alpha, k,\ell})_{(\alpha, k,\ell) \in S} \in \left(\C^*\right)^{S}$ and a family $(f_{\alpha, k,\ell}(z))_{(\alpha, k,\ell) \in S}$ of nonzero $G$-functions. We consider $$f(z)=\sum\limits_{(\alpha, k, \ell) \in S}^{} c_{\alpha, k,\ell} z^{\alpha} \log(z)^k f_{\alpha, k,\ell}(z)$$ a Nilsson-Gevrey series of arithmetic and holonomic type of order $0$. Then $f(z)$ is solution of a nonzero linear differential equation with coefficients in $\Qbar(z)$ and the minimal operator $L$ of $f(z)$ over $\Qbar(z)$ is a $G$-operator.
\end{Th}
 In this paper, we consider only Nilsson-Gevrey series of order $0$, so that we shall not write "order $0$" anymore. 
 
 We aim to provide a more precise result: is it possible to find a quantitative relation between the size $\sigma(L)$ of the minimal operator of $f(z)$ over $\Qbar(z)$ and the sizes of the $G$-functions $f_{\alpha, k,\ell}(z)$?

This question will be studied and given a positive answer in Section \ref{sec:chudnovskynilssongevrey}. The main result of this paper -- and the answer to the problem above -- is the following theorem:

\begin{Th}
\label{th:chudnovskynilssongevreyandre}
We keep the notations of Theorem \ref{th:chudnovskynilssongevreyandreI}. For every $(\alpha, k, \ell) \in S$, let $L_{\alpha, k,\ell} \neq 0$ denote a minimal operator of the $G$-function $f_{\alpha, k,\ell}(z)$. We set $\kappa$ the maximum of  the integers $k$ such that $(\alpha, k, \ell) \in S$ for some $(\alpha, \ell) \in \Q \times \N$ and $A=\{ \alpha \in \Q : \exists (k, \ell) \in \N^2, (\alpha, k, \ell) \in S \}$. Then we have \begin{multline}\label{eq:formulethprinc}
    \quad \sigma(L) \leqslant \max\left(1+\log(\kappa+2), \; 2\big(1+\log(\kappa+2)\big)\log\big(\max\limits_{\alpha \in A} \den(\alpha)\big),\right.\\  \left.\max\limits_{(\alpha, k, \ell) \in S}\big((1+\log(k+2))\sigma(L_{\alpha, k,\ell})\big) \right).\quad
\end{multline}
\end{Th}

It turns out that the size of a differential system is invariant by equivalence of differential systems. Thus, André was able to generalize the notion of size to differential modules and to study the behavior of the size with respect to the usual algebraic operations on differential modules: submodule, quotient, direct sum, tensor product, etc.

We will follow this point of view to solve the problem above and devote Section \ref{sec:sizediffmodule} to some useful reminders on the differential modules and to a synthesis of the results of André on the size of differential modules.

In Section \ref{sec:application}, we will finally give two applications of these results. The first one consists, given two $G$-operators $L_1$ and $L_2$, in finding an explicit inequality between the size of $L_1 L_2$ and the sizes of $L_1$ and $L_2$ (Theorem  \ref{th:tailleproduitgop} in \S \ref{subsec:sizeproductGop}). Then, using \ref{th:chudnovskynilssongevreyandre} and Theorem \ref{th:tailleproduitgop}, we give in \S \ref{subsec:diophproblem} an application consisting in the evaluation of a constant appearing in a Diophantine problem, which is studied in another paper \cite{Lepetit2}.

\medskip

\textbf{Acknowledgements:} I thank T. Rivoal for carefully reading this paper and for his useful comments and remarks that improved it substantially. I also thank the anonymous referee for her/his attentive reading of the manuscript and for her/his relevant comments. 

\medskip

\textbf{Data availability statement:} Data sharing not applicable to this article as no datasets were generated or analysed during the current study.

\textbf{Conflict of interest statement:} On behalf of all authors, the corresponding author states that there is no conflict of interest.

\section{Size of a differential module} \label{sec:sizediffmodule}

\subsection{Reminder on differential modules}

 In this part, we will essentially synthetise the presentation of the theory of differential modules contained in \cite[Chapter 2]{Singer}, that will be useful in the rest of this paper. More details and proofs about differential modules can be found in the above-quoted book.

\begin{defi}
Let $(k, \delta)$ be a differential field. The ring of differential operators over $k$ is denoted by $k[\partial]$ and defined as the set of noncommutative polynomials in $\partial$ with the compatibility condition $$\forall a \in k,  \quad \partial a=a \partial + \delta(a).$$
\end{defi}

This is a left (and right) euclidean ring for the euclidean function $$\ord : L= \sum\limits_{k=0}^{n} a_k \partial^k \mapsto n \quad\quad (\text{where} \;\; a_n \neq 0).$$

\begin{defi}
Let $(k, \delta)$ be a differential field. A differential module $\mathcal{M}$ over $k$ is a $k$-vector space of finite dimension which is a left $k[\partial]$-module, or equivalently a $k$-vector space of finite dimension endowed with a map $\partial : \mathcal{M} \rightarrow \mathcal{M}$ such that $$\forall \lambda \in k, \;\; \forall a, b \in \mathcal{M}, \quad \partial(a+b)=\partial(a)+\partial(b) \quad \text{and} \quad \partial(\lambda a)=\lambda\partial(a)+\delta(\lambda)a.$$
\end{defi}

Let $\mathcal{M}$ be a $k$-differential module of dimension $n$ endowed with a $k$-basis $(e_1, \dots, e_n)$. Then there exists a matrix $A=(a_{i,j})_{i,j} \in M_n(k)$ such that $\forall 1 \leqslant i \leqslant n, \; \partial e_i=-\sum\limits_{j=1}^n a_{j,i} e_j$. A direct computation shows that if $u=\sum\limits_{i=1}^n u_i e_i \in \mathcal{M}$, then $\partial u=0 \iff u'=Au, u={}^t (u_1, \dots, u_n)$.

\begin{rqu}
If we choose another $k$-basis $(f_1, \dots, f_n)$ of $\mathcal{M}$ such that $\partial f_i=-\sum\limits_{j=1}^n b_{j,i} f_j$, then there exists $P \in \GL_n(k)$ such that $B=P[A]=P^{-1}P'+P^{-1}AP$. 

The differential systems defined by $A$ and $B$ are then said to be \emph{equivalent} over $k$. This is equivalent to the fact that $y \mapsto Py$ is a one-to-one correspondance between the solution sets of $y'=Ay$ and $y'=By$.

Therefore, the matrix $A$ obtained by the construction above doesn't depend on the choice of basis we make, up to equivalence of differential systems.
\end{rqu}

We can also do the converse operation:

\begin{defi}
For $A=(a_{ij})_{i,j} \in M_n(k)$, we define the differential module $\mathcal{M}_A$ associated with the differential system $y'=Ay$ by $\mathcal{M}_{A}=k^n$ endowed with the derivation $\partial$ such that $\forall i \in \{ 1, \dots, n \}$, $\partial e_i=-\sum\limits_{j=1}^n a_{ji} e_j$, where $(e_1, \dots, e_n)$ is the canonical basis of $k^n$.
\end{defi}

We associate with a differential operator $L=\partial^n+\sum\limits_{k=0}^n a_k \partial^k \in k[\partial]$ the differential module $\mathcal{M}_L := \mathcal{M}_{A_L}$, where $A_L$ is the companion matrix of $L$: $$A_L=\begin{pmatrix}
0 & 1 &  & (0) \\ 
 & \ddots & \ddots &  \\ 
(0) &  & 0 & 1 \\ 
-a_0 & \dots &  & -a_{n-1}
\end{pmatrix}.$$

\begin{defi}
A \emph{morphism of differential modules} is a $k$-linear map $\varphi : (\mathcal{M}, \partial_{\mathcal{M}}) \rightarrow (\mathcal{N}, \partial_{\mathcal{N}})$ such that $\partial_{\mathcal{N}} \circ \varphi=\varphi \circ \partial_{\mathcal{M}}$.
\end{defi}

\begin{prop}\label{prop:isomodule=equiv}
Let $A, B \in M_n(k)$. Then there is an isomorphism of differential modules betweeen $\mathcal{M}_A$ and $\mathcal{M}_B$ if and only if $A$ and $B$ are equivalent over $k$.
\end{prop}

Thus, the classes of equivalences of differential systems over $k$ classify the $k$-differential modules up to isomorphism.

\bigskip

We can perform usual algebraic constructions with the differential modules:
\begin{itemize}
    \item A \emph{differential submodule} of $\mathcal{M}$ is a left $k[\partial]$-submodule of $\mathcal{M}$.
    \item If $\mathcal{N}$ is a differential submodule of $\mathcal{M}$, then the \emph{quotient differential module} $\mathcal{M}/\mathcal{N}$ is endowed with the quotient derivation $\overline{\partial} : \overline{m} \mapsto \overline{\partial(m)}$.
    \item The \emph{cartesian product} (or direct sum) of $(\mathcal{M}_1, \partial_1)$ and $(\mathcal{M}_2, \partial_2)$ can be endowed with the derivation $\partial : (m_1, m_2) \mapsto (\partial_1(m_1), \partial_2(m_2))$, which makes it a differential module.
    \item The \emph{tensor product} $\mathcal{M}_1 \otimes_k \mathcal{M}_2$ is a differential module when endowed with $\partial : m_1 \otimes m_2 \mapsto m_1 \otimes \partial_2(m_2) + \partial_1(m_1) \otimes m_2$ (extended on the whole vector space with the additivity of $\partial$). Note that the tensor product cannot be defined over the noncommutative ring $k[\partial]$, since $\mathcal{M}_1$ and $\mathcal{M}_2$ are only left $k[\partial]$-modules. We will therefore often denote $\mathcal{M}_1 \otimes \mathcal{M}_2$ without precising the base field.
    \item The derivation $\partial^* : \varphi \mapsto \left(m \mapsto \varphi(\partial(m))-\delta(\varphi(m))\right)$ makes the \emph{dual} $\mathcal{M}^*=\mathrm{Hom}_k(M,k)$ a differential module 
    \item More generally, the preceding constructions enable us to define a structure of differential module on the \emph{space of morphisms} $\mathrm{Hom}_k(\mathcal{M}_1, \mathcal{M}_2) \simeq \mathcal{M}_1^* \otimes_k \mathcal{M}_2$.
\end{itemize}

\begin{rqu}
Actually, for any $A \in M_n(k)$, we have $\mathcal{M}_A^* \simeq \mathcal{M}_{A^*}$, where $A^*=-{}^t A$. If $X$ is a fundamental matrix of solutions of the system $y'=Ay$, then $Y:= {}^t X^{-1}$ satisfies $Y'=-{}^t A Y$. 
\end{rqu}

\begin{lem}
There is an isomorphism of differential modules $\mathcal{M}_L \simeq \left(\Qbar(z)[\partial]/\Qbar(z)[\partial]L\right)^*$.
\end{lem}

\begin{defi}
The \emph{adjoint operator} of $L \in \Qbar(z)[\partial]$ is the operator $L^*$ such that $\mathcal{M}_{L^{*}} \simeq \Qbar(z)[\partial]/\Qbar(z)[\partial]L$.
\end{defi}

We can compute directly $L^*$: if $L=\partial^n+\sum\limits_{k=0}^{n-1} a_k \partial^k$, then $L^*=(-\partial)^n+\sum\limits_{k=0}^{n-1} (-\partial)^k a_k$. The adjoint operator has useful properties, for instance that $L^{**}=L$ and $\left(L_1L_2\right)^*=L_2^* L_1^*$.

\bigskip

We finally mention that the \emph{cyclic vector Theorem}, implies that any differential system $y'=Ay$ is equivalent to some system $y'=A_L y$, $L \in k[\partial]$. This yields \emph{in fine} an equivalence between differential systems and differential operators (see \cite[pp. 42--43]{Singer}).

\subsection{Size and operations on differential modules} \label{subsec:Sizeandops}

Let $\K$ be a number field, and $G \in M_n(\K(z))$. In this subsection, following André's notations, we introduce a quantity $\sigma(G)$ encoding an arithmetic condition of moderate growth on a differential system, called \emph{Galochkin condition}. 

For $s \in \N$, we define $G_s$ as the matrix such that, if $y$ is a vector satisfying $y'=Gy$, then $y^{(s)}=G_s y$. In particular, $G_0$ is the identity matrix. The matrices $G_s$ satisfy the recurrence relation $$\forall s \in \N, \quad G_{s+1}=G_sG+G'_s,$$ where $G'_s$ is the derivative of the matrix $G_s$.

\begin{defi}[Galochkin, \cite{Galochkin74}] \label{def:galochkin}
Let $T(z) \in \K(z)$ be such that $T(z)G(z) \in M_n(\K[z])$. For $s \in \N$, we consider $q_s \geqslant 1$ the least common denominator of all coefficients of the entries of the matrices $T(z)^m \gf{G_m(z)}{m!}$, when $m \in \{ 1, \dots, s \}$. The system $y'=Gy$ is said to \emph{satisfy the Galochkin condition} if $$\exists C>0: \; \forall s \in \N, \quad  q_s \leqslant C^{s+1}.$$
\end{defi}

Chudnovsky's Theorem (cf \cite{Chudnovsky}) states that if $G$ is the companion matrix of the minimal non-zero differential operator $L$ associated to a $G$-function, then the system $y'=Gy$ satisfies the Galochkin condition. That is why we say that $L$ is a \emph{$G$-operator} (see \cite[pp. 717--719]{AndregevreyI} for a review of the properties of $G$-operators). Following \cite[chapter VII]{Dwork}, we are now going to rephrase this condition in $p$-adic terms.

For every prime ideal $\mathfrak{p}$ of $\Oal_{\K}$, we define $|\cdot|_{\mathfrak{p}}$ as the $\mathfrak{p}$-adic absolute value on $\K$, with the choice of normalisations given in \cite[p.223]{Dwork}.

We recall that the \emph{Gauss absolute value} associated with $|\cdot|_{\mathfrak{p}}$ is the non-archimedean absolute value $$\fonction{|\cdot|_{\mathfrak{p},\mathrm{Gauss}}}{\K(z)}{\R}{\gf{\sum\limits_{i=0}^{N} a_i z^i}{\sum\limits_{j=0}^{M} b_j z^j}}{\gf{\max\limits_{0 \leqslant i \leqslant N} |a_i|_{\mathfrak{p}}}{\max\limits_{0 \leqslant j \leqslant M} |b_j|_{\mathfrak{p}}}.}$$ 

The absolute value $|\cdot|_{\mathfrak{p},\mathrm{Gauss}}$ naturally induces a norm on $M_{n,m}(\K(z))$, defined for all $H=(h_{i,j})_{i,j} \in M_{n,m}(\K(z))$ as $\|H\|_{\mathfrak{p},\mathrm{Gauss}}=\max_{i,j} |h_{i,j}|_{\mathfrak{p},\mathrm{Gauss}}$. It is called the \emph{Gauss norm}. If $n=m$, $\|\cdot\|_{\mathfrak{p},\mathrm{Gauss}}$ is a norm of algebra on $M_n(\K(z))$. We now use this notation to define the notion of size of a matrix:

\begin{defi}[\cite{Dwork}, p. 227]\label{def:taillematrice}
Let $G \in M_{n}(\K(z))$. The \emph{size} of $G$ is $$\sigma(G):=\limsup\limits_{s \rightarrow + \infty} \gf{1}{s} \sum\limits_{\mathfrak{p} \in \Spec(\Oal_\K)} h(s, \mathfrak{p})$$ where $$ \forall s \in \N, \quad h(s, \mathfrak{p})= \sup_{m \leqslant s} \log^{+} \left\| \gf{G_m}{m!} \right\|_{\mathfrak{p}, \mathrm{Gauss}},$$  with $\log^{+} : x \mapsto \log\left(\max(1,x)\right)$. The \emph{size} of $Y=\sum\limits_{m=0}^{\infty} Y_m z^m, Y_m  \in M_{p,q}(\K)$ is $$\sigma(Y):=\limsup\limits_{s \rightarrow + \infty} \gf{1}{s} \sum\limits_{\mathfrak{p} \in \Spec(\Oal_\K)} \sup_{m \leqslant s} \log^{+} \left\| Y_m \right\|_{\mathfrak{p}}. $$
\end{defi}

The relation between Galochkin condition and size is given by the following result. 

\begin{prop} \label{prop:lienq'ssigmaG} 
\begin{enumerateth}
 \item With the notations of Definitions \ref{def:galochkin} and \ref{def:taillematrice}, we have $$ \sigma(G)+\gf{h^{-}(T)}{[\K : \Q]} \leqslant \limsup\limits_{s \rightarrow + \infty} \gf{1}{s} \log(q_s)  \leqslant [\K:\Q] \sigma(G) + h^{+}(T)\;,$$
 where $$h^{-}(T)=\sum\limits_{\mathfrak{p} \in \Spec(\Oal_\K)} \log\left(\min\left(1, |T|_{\mathfrak{p}, \mathrm{Gauss}}\right)\right) \quad \mathrm{and} \quad h^{+}(T)=\sum\limits_{\mathfrak{p} \in \Spec(\Oal_\K)} \log^{+} |T|_{\mathfrak{p}, \mathrm{Gauss}}$$ \cite[Proposition 5, p. 16]{Lepetit2}.
 \item   In particular, the differential system $y'=Gy$ satisfies the Galochkin condition if and only if $\sigma(G) < + \infty$ \cite[p. 228]{Dwork}.
\end{enumerateth}

\end{prop}

When $G \in M_{n}(\Q(z))$ and $T(z) \in \Z[z]$ has at least one coefficient equal to $1$, we have the equality $$\sigma(G)=\limsup\limits_{s \rightarrow + \infty} \gf{1}{s} \log(q_s).$$

The size does not depend on the number field $\K$ we consider, and it is invariant under equivalence of differential systems, as this lemma shows. It is due to André \cite[p. 71]{Andre}. We give a more detailed proof than in \cite{Andre} for the sake of completeness.

\begin{lem}[\cite{Andre}, p. 71]\label{lem:sizeinvariantequivsystdiff}
Let $\K$ a number field, $G \in M_n(\K(z))$ and $P \in \GL_n(\K(z))$, let $H=P[G]=PGP^{-1}+P'P^{-1}$ be a matrix defining a differential system $y'=Hy$ which is equivalent to $y'=Gy$ over $\Qbar(z)$. Then $\sigma(G)=\sigma(H)$.
\end{lem}

\begin{dem}
We have$$\forall s \in \N, \; \; H_s=\sum_{m=0}^s \dbinom{s}{m} P^{(s-m)} G_m P^{-1}.$$ 

Indeed, let $s \in \N$ and $y$ such that $y'=Gy$. Then $y^{(s)}=G_s y$ and, as $H=P[G]$, $(Py)^{(s)}=H_s Py$. But $$H_s Py= (Py)^{(s)}=\sum_{k=0}^s \binom{s}{k} P^{(s-k)} y^{(k)}=\sum_{k=0}^s \binom{s}{k} P^{(s-k)} G_k y.$$ Since the equality holds for every solution $y$ of $y'=Gy$, it follows that $$H_s=\sum_{m=0}^s \binom{s}{m} P^{(s-m)} H_m P^{-1}.$$

Therefore, for $\mathfrak{p} \in \Spec(\Oal_{\K})$, \begin{align*}\left\|\gf{H_s}{s!}\right\|_{\mathfrak{p},\mathrm{Gauss}} &\leqslant \max_{0 \leqslant m \leqslant s} \left\|\gf{P^{(s-m)}}{(s-m)!}\right\|_{\mathfrak{p},\mathrm{Gauss}} \cdot  \|P^{-1}\|_{\mathfrak{p},\mathrm{Gauss}} \cdot \left\|\gf{G_m}{m!}\right\|_{\mathfrak{p}, \mathrm{Gauss}} \\
&\leqslant \|P\|_{\mathfrak{p},\mathrm{Gauss}} \cdot \|P^{-1}\|_{\mathfrak{p},\mathrm{Gauss}} \cdot \max_{0 \leqslant m \leqslant s} \left\|\gf{G_m}{m!}\right\|_{\mathfrak{p}, \mathrm{Gauss}} \end{align*} because $\left\|\gf{P^{(s-m)}}{(s-m)!}\right\|_{\mathfrak{p},\mathrm{Gauss}} \leqslant \|P\|_{\mathfrak{p}, \mathrm{Gauss}}$ (cf \cite[p. 118]{Dwork}). Denoting $C(\mathfrak{p})=\|P\|_{\mathfrak{p},\mathrm{Gauss}} \cdot \|P^{-1}\|_{\mathfrak{p},\mathrm{Gauss}}$, we then have, for all integers $s$, $$\max_{0 \leqslant m \leqslant s} \left\|\gf{H_s}{s!}\right\|_{\mathfrak{p},\mathrm{Gauss}} \leqslant C(\mathfrak{p}) \max_{0 \leqslant m \leqslant s} \left\|\gf{G_m}{m!}\right\|_{\mathfrak{p}, \mathrm{Gauss}}.$$ Hence $$\sigma(H) \leqslant \sigma(G)+\limsup_{s \rightarrow + \infty} \gf{1}{s} \sum_{\mathfrak{p} \in \Spec(\Oal_{\K})} \log^{+} C(\mathfrak{p})= \sigma(G).$$ Symmetrically, we have $\sigma(G) \leqslant \sigma(H)$, as $G=P^{-1}[H]$. 
\end{dem}

Lemma \ref{lem:sizeinvariantequivsystdiff} and Proposition \ref{prop:isomodule=equiv} imply that we can define without ambiguity the \emph{size} of a differential module $\mathcal{M}$ by $\sigma(\mathcal{M}) := \sigma(A)$, because every differential module $\mathcal{M}$ can be associated with a unique class of equivalence of differential systems $[A]$.

\begin{defi}[\cite{Andre}, pp. 74--76]
Let $L \in \Qbar(z)[\partial]$. We set $\sigma(L)=\sigma(\mathcal{M}_{L})$. If $\sigma(L) < +\infty$, the operator $L$ is said to be a \emph{$G$-operator}.
\end{defi}

This terminology is justified by the fact that every solution of the equation $L(y(z))=0$ around a non singular point of the $G$-operator $L$ is a $G$-function.

\bigskip

The following result is key to this paper. It is due to André. It proves that the notion of size is compatible with most of the usual operations on differential modules.

Recall that for $n \in \N^*$, the $n$-th symmetric power of the differential module $\mathcal{M}$, $\mathrm{Sym}^n_{\Qbar(z)}(\mathcal{M})$, is the quotient of $\mathcal{M}^{\otimes n}$ by the submodule generated over $\Qbar(z)$ by the $m_1 \otimes \dots \otimes m_n - m_{\sigma(1)} \otimes \dots \otimes m_{\sigma(n)}$, when $\sigma$ is any permutation of $\{1, \dots, n \}$.

\begin{prop}[\cite{Andre}, p. 72] \label{prop:tailleetopmodules}
Let $\mathcal{M}_1, \mathcal{M}_2, \mathcal{M}_3$ be differential modules over $\Qbar(z)$. Then 
\begin{enumerateth}
\item If $\mathcal{M}_2$ is a differential submodule of $\mathcal{M}_1$, then $\sigma(\mathcal{M}_2) \leqslant \sigma(\mathcal{M}_1)$ and $\sigma\left(\mathcal{M}_1 \middle/\mathcal{M}_2\right) \leqslant \sigma(\mathcal{M}_1)$.
\item $\sigma(\mathcal{M}_1 \times \mathcal{M}_2)=\max\left(\sigma(\mathcal{M}_1), \sigma(\mathcal{M}_2)\right) \leqslant \sigma(\mathcal{M}_1) + \sigma(\mathcal{M}_2)$.
\item $\sigma\left(\Sym^N_{\Qbar(z)}(\mathcal{M}_1)\right) \leqslant \left(1+\log(N)\right) \sigma(\mathcal{M}_1)$.
\item $\sigma(\mathcal{M}_1^*) \leqslant \sigma(\mathcal{M}_1)\big(1+\log(\mu_1-1)\big)$, where $\mu_1=\dim_{\Qbar(z)}\mathcal{M}_1$.
\item If there is an exact sequence $0 \rightarrow \mathcal{M}_2 \rightarrow \mathcal{M}_1 \rightarrow \mathcal{M}_3 \rightarrow 0$, then we have \begin{equation} \label{eq:seqexacteeq1}\sigma(\mathcal{M}_1) \leqslant 1+2\sigma(\mathcal{M}_2 \times \mathcal{M}_3\times (\mathcal{M}_3)^*) \leqslant 1+2\max\left(\sigma(\mathcal{M}_2), \sigma(\mathcal{M}_3\right)+\mu_3-1)\end{equation} with $\mu_3=\dim_{\Qbar(z)} \mathcal{M}_3$.
\end{enumerateth}
\end{prop}
 
In \textbf{d)}, we have the following alternative bound that we will show in the proof of Proposition \ref{prop:tailleetopmodules}: \begin{equation}\label{eq:sizeadjalt}\sigma(\mathcal{M}_1)-\mu_1+1 \leqslant \sigma(\mathcal{M}_1^*) \leqslant \sigma(\mathcal{M}_1)+\mu_1-1.\end{equation}

In \textbf{e)}, we can actually deduce from André's proof the following slightly more precise result: \begin{equation}\label{eq:seqexacteeq2}\sigma(\mathcal{M}_1) \leqslant 1+\gf{11}{6}\sigma\left(\mathcal{M}_2 \times \mathcal{M}_3\times (\mathcal{M}_3)^*\right).\end{equation} 
Alternatively (see \eqref{eq:sizeadjalt} above and \eqref{eq:result1tailleopmod} below), we have \begin{equation} \label{eq:seqexacteeq3}\sigma(\mathcal{M}_1) \leqslant 1+\sigma(\mathcal{M}_2)+\sigma(\mathcal{M}_3)+\sigma(\mathcal{M}_3^*) \leqslant\mu_3+\sigma(\mathcal{M}_2)+2\sigma(\mathcal{M}_3).\end{equation}

Note that \eqref{eq:seqexacteeq2} is not always a better estimate than \eqref{eq:seqexacteeq3}, as the example with $\sigma(\widetilde{L}_{\beta})$ in Subsection \ref{subsec:diophproblem} shows.

The proof of \textbf{a)} and \textbf{e)} relies essentially on the following lemma, stated without proof in \cite[p. 72]{Andre}, of which we give a proof for the reader's convenience:
\begin{lem}[\cite{Andre}, p. 72] \label{lem:formematriceseqexacte}
Let $\mathcal{M}_1, \mathcal{M}_2, \mathcal{M}_3$ be differential modules over $\Qbar(z)$. If there is an exact sequence $0 \rightarrow \mathcal{M}_2 \rightarrow \mathcal{M}_1 \rightarrow \mathcal{M}_3 \rightarrow 0$, then, with a suitable choice of basis, $\mathcal{M}_1$ is associated with a differential system of the form $\partial y=Gy$, where  $G=\begin{pmatrix} G^{(2)} & G^{(0)} \\ 0 & G^{(3)} \end{pmatrix}$ and the systems $\partial y=G^{(2)} y$ and $\partial y=G^{(3)} y$ are respectively associated with $\mathcal{M}_2$ and $\mathcal{M}_3$.
\end{lem}
\begin{dem}
We read on the exact sequence that $\mathcal{M}_3 \simeq  \left.\mathcal{M}_1 \middle/ \mathcal{M}_2\right.$.

Linear algebra convinces us that, if $(e_1, \dots, e_n)$ is a $k$-basis of $\mathcal{M}_2$ and $(\overline{f_1}, \dots, \overline{f_m})$ is a $k$-basis of $\mathcal{M}_1/\mathcal{M}_2$, then $\mathcal{B}=(e_1, \dots, e_n, f_1, \dots, f_m)$ is a basis of $\mathcal{M}_1$.

Thus, if $\forall 1 \leqslant i \leqslant m, \;\; \partial \overline{f_i}=-\sum\limits_{j=1}^m g^{(3)}_{j,i} \overline{f_j}$, we have $\partial f_i+\sum\limits_{j=1}^n g^{(3)}_{j,i} f_j \in \mathcal{M}_2$ so that there exists a matrix $G^{(0)}=(g^{(0)}_{i,j}) \in M_{m,n}(k)$ such that $$\forall 1 \leqslant i \leqslant m, \quad \partial f_i=-\sum\limits_{j=1}^m g^{(3)}_{j,i} f_j- \sum\limits_{j=1}^n g^{(0)}_{j,i} e_j.$$ Furthermore, we can find a matrix $G^{(2)} \in M_n(k)$ associated with $\mathcal{M}_2$ such that $$\forall 1 \leqslant i \leqslant n, \quad \partial e_i=-\sum\limits_{j=1}^m g^{(2)}_{j,i} e_j.$$ Hence, the matrix $G$ such that $\partial {}^t \mathcal{B}=-{}^t G {}^t \mathcal{B}$ is of the form $G=\begin{pmatrix} G^{(2)} & G^{(0)} \\ 0 & G^{(3)} \end{pmatrix}$.
\end{dem}

\begin{dem}[of Proposition \ref{prop:tailleetopmodules}]

We give this proof for the reader's convenience, providing moreover some details that were not written by André.

\medskip

\textbf{a)} is a direct consequence of Lemma \ref{lem:formematriceseqexacte}, and \textbf{b)} follows straightforwardly from the definition. For the point \textbf{c)}, see \cite[p. 72]{Andre} (and \cite[pp. 17--18]{Andre} for the key argument). 

The inequality \textbf{d)} is a consequence of Katz's Theorem on the exponents of a $G$-operator (see \cite[p. 80]{Andre}). Let us explain how to obtain the alternative bound \eqref{eq:sizeadjalt} on $\sigma(\mathcal{M}_1^*)$.

Following the notations of \cite[p. 226]{Dwork}, we denote by $\rho(\mathcal{M}_1)$ the global inverse radius of $\mathcal{M}_1$. Then the André-Bombieri Theorem \cite[p. 74]{Andre} yields $$\rho(\mathcal{M}_1) \leqslant \sigma(\mathcal{M}_1) \leqslant \rho(\mathcal{M}_1)+\mu_1-1,$$ where $\mu_1=\dim_{\Qbar(z)} \mathcal{M}_1$. Moreover, \cite[Lemma 2 p. 72]{Andre} states that $\rho(\mathcal{M}_1^*)=\rho(\mathcal{M}_1)$. Thus, the application of André-Bombieri Theorem to the adjoint module $\mathcal{M}_1^*$ yields $\rho(\mathcal{M}_1) \leqslant \sigma(\mathcal{M}_1^*) \leqslant \rho(\mathcal{M}_1)+\mu_1-1$, hence $$\sigma(\mathcal{M}_1)-\mu_1+1 \leqslant \sigma(\mathcal{M}_1^*) \leqslant \sigma(\mathcal{M}_1)+\mu_1-1.$$

\bigskip

We now prooced to the proof of \textbf{e)}. André gave a proof of this assertion without providing all the details. We present them here. 

By Lemma \ref{lem:formematriceseqexacte}, we can find suitable bases of $\mathcal{M}_1$, $\mathcal{M}_2$, $\mathcal{M}_3$ such that, in these bases, $\mathcal{M}_1$ (resp. $\mathcal{M}_2$, $\mathcal{M}_3$) represents the system $\partial y=Gy$ (resp. $\partial y=G^{(2)} y$ and $\partial y=G^{(3)} y$) where $$G=\begin{pmatrix} G^{(2)} & G^{(0)} \\ 0 & G^{(3)} \end{pmatrix}.$$

Let $\K$ be a number field such that $G \in M_n(\K(z))$.
We fix $\mathfrak{p} \in \Spec(\Oal_{\K})$. We consider $t_{\mathfrak{p}}$ a free variable on $\K$ called the \emph{generic point}. We can then build a complete and algebraically closed extension of $\left(\K(t_{\mathfrak{p}}), | \cdot |_{\mathfrak{p}, \mathrm{Gauss}}\right)$. Concretely, $\Omega_{\mathfrak{p}}$ is the completion of the algebraic closure of the completion of $\K(t_{\mathfrak{p}})$ endowed with the Gauss valuation $| \cdot |_{\mathfrak{p}, \mathrm{Gauss}}$ (see \cite[p. 93]{Dwork}).

We set $$X^{(2)}_{t_{\mathfrak{p}}}(z)=\sum\limits_{s=0}^{\infty} \gf{G^{(2)}_s(t_{\mathfrak{p}})}{s!} (z-t_{\mathfrak{p}})^s \in \GL_n(\Omega_{\mathfrak{p}}\llbracket z-t_{\mathfrak{p}}\rrbracket)$$ (resp. $X^{(3)}_{t_{\mathfrak{p}}} \in \GL_m(\Omega_{\mathfrak{p}}\llbracket z-t_{\mathfrak{p}}\rrbracket)$) a fundamental matrix of solutions of $\partial y=G^{(2)} y$ (resp. $\partial y=G^{(3)}y$) at the generic point $t_{\mathfrak{p}}$ such that $X^{(2)}_{t_{\mathfrak{p}}}(t_{\mathfrak{p}})=I_n$ (resp. $X^{(3)}_{t_{\mathfrak{p}}}(t_{\mathfrak{p}})=I_m$). Consider $X^{(0)}_{t_{\mathfrak{p}}} \in M_{n,m}(\Omega_{\mathfrak{p}}\llbracket z-t_{\mathfrak{p}}\rrbracket)$ a solution of \begin{equation} \label{eq:EDX0}\partial X_{t_{\mathfrak{p}}}^{(0)}=G^{(3)}X_{t_{\mathfrak{p}}}^{(0)}+G^{(0)} X_{t_{\mathfrak{p}}}^{(2)} \end{equation} such that $X^{(0)}_{t_{\mathfrak{p}}}(t_{\mathfrak{p}})=0$.

Then $X_{t_{\mathfrak{p}}}=\begin{pmatrix} X_{t_{\mathfrak{p}}}^{(2)} & X_{t_{\mathfrak{p}}}^{(0)} \\ 0 & X_{t_{\mathfrak{p}}}^{(3)} \end{pmatrix}$ is a fundamental matrix of solutions of $\partial y=Gy$ such that $X_{t_{\mathfrak{p}}}(t_{\mathfrak{p}})=I_{n+m}$. Thus we have $\forall s \in \N,\; \partial^s(X_{t_{\mathfrak{p}}})(t_{\mathfrak{p}})=G_s(t_{\mathfrak{p}})$. Since we know that $\partial^s(X^{(2)}_{t_{\mathfrak{p}}})(t_{\mathfrak{p}})=G^{(2)}_s(t_{\mathfrak{p}})$ (and likewise for $X_{t_{\mathfrak{p}}}^{(3)}$), it  suffices to estimate the Gauss norm of $\partial^s(X^{(0)}_{t_{\mathfrak{p}}})(t_{\mathfrak{p}})$ to obtain an estimate on the size of $G$.

For simplicity, we will omit the index $t_{\mathfrak{p}}$ in what follows.

It follows from \eqref{eq:EDX0} that \begin{align*}
    \partial(X^{(3)^{-1}} X^{(0)})&=X^{(3)^{-1}} \partial(X^{(0)})-X^{(3)^{-1}} \partial(X^{(3)}) X^{(3)^{-1}} X^{(0)} \\ &=X^{(3)^{-1}} (G^{(3)}X^{(0)}+G^{(0)} X^{(2)})-X^{(3)^{-1}} G^{(3)} X^{(0)} =X^{(3)^{-1}} G^{(0)} X^{(2)} 
\end{align*} so that, for $\ell \in \N$, by the Leibniz formula, \begin{align*}
    \gf{\partial^\ell X^{(0)}}{\ell!} &=\sum\limits_{i=0}^\ell  \gf{\partial^{\ell-i}(X^{(3)})}{(\ell-i)!} \gf{\partial^i (X^{(3)^{-1}} X^{(0)})}{i!}=\gf{\partial^\ell (X^{(3)})}{\ell!} X^{(3)^{-1}} X^{(0)}+ \sum\limits_{i=1}^\ell  \gf{\partial^{\ell-i}(X^{(3)})}{(\ell-i)!} \gf{\partial^{i-1} (X^{(3)^{-1}} G^{(0)} X^{(2)})}{i!} \\
    &= \gf{partial^\ell (X^{(3)})}{\ell!} X^{(3)^{-1}} X^{(0)}+ \sum\limits_{i=0}^\ell  \gf{\partial^{\ell-i}(X^{(3)})}{(\ell-i)!} \gf{1}{i} \sum_{i_1+i_2+i_3=i-1} \gf{\partial^{i_1}(X^{(3)^{-1}})}{i_1!} \gf{\partial^{i_2}(G^{(0)})}{i_2!} \gf{\partial^{i_3}(X^{(2)})}{i_3!}.
\end{align*}

We now evaluate this equality at the generic point $t_{\mathfrak{p}}$. Since $X^{(0)}(t_{\mathfrak{p}})=0$, we get
$$\gf{\partial^\ell X^{(0)}}{\ell!}(t_{\mathfrak{p}})= \sum\limits_{i=0}^\ell  \gf{\partial^{\ell-i}(X^{(3)})(t_{\mathfrak{p}})}{(\ell-i)!} \gf{1}{i} \sum_{i_1+i_2+i_3=i-1} \gf{\partial^{i_1}(X^{(3)^{-1}})(t_{\mathfrak{p}})}{i_1!} \gf{\partial^{i_2}(G^{(0)})(t_{\mathfrak{p}})}{i_2!} \gf{\partial^{i_3}(X^{(2)})(t_{\mathfrak{p}})}{i_3!}.$$

We make the following observations:
\begin{itemize}
    \item Since $G^{0}(t_{\mathfrak{p}})$ has coefficients in $\K(t_{\mathfrak{p}})$, we have $\left\| \gf{\partial^{i_2}(G^{(0)})}{i_2!} \right\|_{\mathfrak{p},\mathrm{Gauss}} \leqslant \|G^{0} \|_{\mathfrak{p}, \mathrm{Gauss}}$ by \cite[p. 94]{Dwork}.
    
    \item Set $\delta_\ell=\mathrm{lcm}(1,2, \dots, \ell)$. Then  $$\sum\limits_{\mathfrak{p} \in \Spec(\Oal_{\K})} \sup\limits_{1 \leqslant i \leqslant \ell} \gf{1}{|i|_\mathfrak{p}}=\gf{1}{[\K:\Q]} \log(\delta_\ell) \leqslant \gf{\ell}{[\K:\Q]}(1+o(1)),$$ by \cite[Lemma 1.1 p. 225]{Dwork}.
\end{itemize}

Hence, since $\| \cdot \|_{\mathfrak{p}, \mathrm{Gauss}}$ is non-archimedean, we have
\begin{multline}\label{eq:inegalitenormesdegauss}
   \log^{+}\left\|\gf{\partial^\ell X^{(0)}}{\ell!}\right\|_{\mathfrak{p}, \mathrm{Gauss}} \leqslant \max_{0\leqslant i \leqslant n \atop i_1+i_2+i_3=i-1} \left( \log^{+}\left\|\gf{\partial^{\ell-i}(X^{(3)})}{(\ell-i)!}\right\|_{\mathfrak{p}, \mathrm{Gauss}}+ \log^{+}\left\|\gf{1}{i} \right\|_{\mathfrak{p}, \mathrm{Gauss}} +\right.\\ \left.\log^{+}\left\|\gf{\partial^{i_1}(X^{(3)^{-1}})}{i_1!}\right\|_{\mathfrak{p}, \mathrm{Gauss}}+ \log^{+}\left\|\gf{\partial^{i_2}(G^{(0)})}{i_2!}\right\|_{\mathfrak{p}, \mathrm{Gauss}}+ \log^{+}\left\|\gf{\partial^{i_3}(X^{(2)})}{i_3!}\right\|_{\mathfrak{p}, \mathrm{Gauss}}\right)
\end{multline}

But  $$\left\|\gf{\partial^{\ell-i}(X^{(3)^{-1}})}{(\ell-i)!}\right\|_{\mathfrak{p}, \mathrm{Gauss}}=\left\|\gf{\left(G^{(3)}\right)_{\ell-i}}{(\ell-i)!}\right\|_{\mathfrak{p}, \mathrm{Gauss}} \quad \text{and} \quad \left\|\gf{\partial^{i_3}(X^{(2)})}{i_3!}\right\|_{\mathfrak{p}, \mathrm{Gauss}}=\left\|\gf{\left(G^{(2)}\right)_{i_3}}{i_3!}\right\|_{\mathfrak{p}, \mathrm{Gauss}}.$$ Moreover we have $$\left\|\gf{\partial^{i_1}(X^{(3)^{-1}})}{i_1!}\right\|_{\mathfrak{p}, \mathrm{Gauss}}=\left\|\gf{\left(-{}^t G^{(3)}\right)_{i_1}}{i_1!}\right\|_{\mathfrak{p}, \mathrm{Gauss}},$$ because ${}^t X^{(3)^{-1}}$ is a fundamental matrix of solutions of the system $y'=-{}^t G^{(3)} y$, associated with the dual module $\mathcal{M}_3^*$. This yields for every $\ell \in \N^*$ and $N \geqslant \ell$ \begin{multline*}
    \log^{+} \left\|\gf{\partial^{\ell}(X^{(0)})(t_{\mathfrak{p}})}{\ell!}\right\|_{\mathfrak{p}, \mathrm{Gauss}} \leqslant h(N,\mathfrak{p},G^{(3)})+h(N,\mathfrak{p},-{}^t G^{(3)})+ h(N,\mathfrak{p},G^{(2)}) \\ + \log^{+} \| G^{(0)} \|_{\mathfrak{p}, \mathrm{Gauss}}+ \log \max_{0 \leqslant i \leqslant N} \gf{1}{|i|_{p}}.
\end{multline*}

We therefore obtain $$\limsup_{N \rightarrow + \infty} \gf{1}{N} \sum\limits_{\mathfrak{p} \in \Spec(\Oal_{\K})} \sup\limits_{n \leqslant N} \log^{+} \left\|\gf{\partial^{\ell}(X^{(0)})(t_{\mathfrak{p}})}{\ell!}\right\|_{\mathfrak{p}, \mathrm{Gauss}} \leqslant \gf{1}{[\K:\Q]}+ \sigma(\mathcal{M}_3)+\sigma(\mathcal{M}_3^*)+\sigma(\mathcal{M}_2)$$ since $\| G^{(0)} \|_{\mathfrak{p}, \mathrm{Gauss}}=1$ for all but a finite number of primes. Furthermore, \begin{align*}\left\|\gf{G_\ell}{\ell!}\right\|_{\mathfrak{p}, \mathrm{Gauss}}&=\left\|\gf{\partial^{\ell}(X)(t_{\mathfrak{p}})}{\ell!}\right\|_{\mathfrak{p}, \mathrm{Gauss}} \\ &=\max\left(\left\|\gf{\partial^{\ell}(X^{(0)})(t_{\mathfrak{p}})}{\ell!}\right\|_{\mathfrak{p}, \mathrm{Gauss}}, \left\|\gf{\partial^{\ell}(X^{(2)})(t_{\mathfrak{p}})}{\ell!}\right\|_{\mathfrak{p}, \mathrm{Gauss}}, \left\|\gf{\partial^{\ell}(X^{(3)})(t_{\mathfrak{p}})}{\ell!}\right\|_{\mathfrak{p}, \mathrm{Gauss}} \right)\end{align*} so finally \begin{equation}\label{eq:result1tailleopmod}\sigma(\mathcal{M}_1) \leqslant 1+\sigma(\mathcal{M}_3)+\sigma(\mathcal{M}_3^*)+\sigma(\mathcal{M}_2).\end{equation}

\bigskip

 We can actually obtain a more refined inequality by using the following argument. Let $\ell \in \N$ and $i_1+ i_2+ i_3 \leqslant \ell-1$. We assume for example $i_1 \geqslant i_2 \geqslant i_3$, the other cases can be treated likewise and lead to the same final conclusion. By considering separately the cases $i_1 \leqslant \gf{\ell-1}{2}$ and $i_1 > \gf{\ell-1}{2}$, we prove that $i_2 \leqslant  \gf{\ell-1}{2}$ and $i_3 \leqslant \gf{\ell-1}{3}$.
 
 Thus we obtain in inequality \eqref{eq:inegalitenormesdegauss} \begin{multline*}
     \max \left(\log^{+}\left\|\gf{\partial^{i_1}(X^{(3)^{-1}})}{i_1!}\right\|_{\mathfrak{p}, \mathrm{Gauss}}+ \log^{+}\left\|\gf{\partial^{i_2}(G^{(0)})}{i_2!}\right\|_{\mathfrak{p}, \mathrm{Gauss}}+ \log^{+}\left\|\gf{\partial^{i_3}(X^{(2)})}{i_3!}\right\|_{\mathfrak{p}, \mathrm{Gauss}}\right) \\ \leqslant h(\ell,\mathfrak{p},G^{(3)})+h\left(\gf{\ell}{2}, \mathfrak{p}, -{}^t G^{(3)}\right)+h\left(\gf{\ell}{3}, \mathfrak{p}, G^{(2)}\right)
 \end{multline*} so that $$\gf{1}{\ell} \sum\limits_{\mathfrak{p} \in \Spec(\Oal_{\K})} h(\ell,\mathfrak{p},G) \leqslant \gf{\log(\delta_\ell)}{\ell}+\gf{1}{\ell} \sum\limits_{\mathfrak{p} \in \Spec(\Oal_{\K})} \log^{+}\| G^{(0)} \|_{\mathfrak{p},\mathrm{Gauss}}+ v_\ell + \gf{1}{\ell} \left\lfloor\gf{\ell}{2}\right\rfloor v_{\left\lfloor\frac{\ell}{2}\right\rfloor}+\gf{1}{\ell} \left\lfloor\gf{\ell}{3}\right\rfloor v_{\left\lfloor\frac{\ell}{3}\right\rfloor} $$ with $v_\ell=\gf{1}{\ell} \sum\limits_{\mathfrak{p} \in \Spec(\Oal_{\K})} \max\left(h(\ell,\mathfrak{p},G^{(3)}), h(\ell,\mathfrak{p},-{}^t G^{(3)}), h(\ell,\mathfrak{p},G^{(2)})\right)$. Taking the superior limit on both sides, we get  $$\sigma(\mathcal{M}_1)\leqslant 1+ \gf{11}{6} \max\left(\sigma(\mathcal{M}_3^*), \sigma(\mathcal{M}_3), \sigma(\mathcal{M}_2)\right),$$ bounded by $1+2 \max\left(\sigma(\mathcal{M}_3^*), \sigma(\mathcal{M}_3), \sigma(\mathcal{M}_2)\right)$ in \cite{Andre}. This is the desired conclusion.\end{dem}

\medskip

Points \textbf{a)}, \textbf{b)} and \textbf{c)} of Proposition \ref{prop:tailleetopmodules} actually imply the following result, which was not stated explicitly by André and is important for our application to Theorem \ref{th:chudnovskynilssongevreyandreI}.

\begin{prop}\label{prop:tailleetproduittensoriel}
Let $\left(\mathcal{M}_i\right)_{1 \leqslant i \leqslant N}$ be a family of differential modules over $\Qbar(z)$. Then $$\sigma(\mathcal{M}_1 \otimes_{\Qbar(z)} \dots \otimes_{\Qbar(z)} \mathcal{M}_N) \leqslant \big(1+\log(N)\big)\max\big(\sigma(\mathcal{M}_1), \dots, \sigma(\mathcal{M}_N)\big).$$
In particular, we have $\sigma\left(\mathcal{M}_1^{\otimes N}\right) \leqslant (1+\log(N))\sigma(\mathcal{M}_1)$.
\end{prop} 

\begin{dem}
We denote by $\overline{\alpha}$ the class of an element $\alpha \in \left(\mathcal{M}_1 \times \dots \times \mathcal{M}_N \right)^{\otimes N}$ in the quotient space $\Sym^N\left(\mathcal{M}_1 \times \dots \times \mathcal{M}_N \right)$. There is an injective morphism of differential modules $$\fonction{\varphi}{\mathcal{M}_1 \otimes \dots \otimes \mathcal{M}_N}{\Sym^N\left(\mathcal{M}_1 \times \dots \times \mathcal{M}_N \right)}{m_1 \otimes \dots \otimes m_N}{\overline{(m_1, 0, \dots,0) \otimes \dots \otimes (0, \dots, 0, m_N)}.}$$ Indeed, for all $i$, let $\mathcal{B}_i=(e^{(i)}_1, \dots, e^{(i)}_{n_i})$ be a $\Qbar(z)$-basis  of $\mathcal{M}_i$. A $\Qbar(z)$-basis of $\mathcal{M}_1 \otimes \dots \otimes \mathcal{M}_N$ is the family $\mathcal{F}$ whose elements are the $e^{(1)}_{i_1} \otimes e^{(2)}_{i_2} \otimes \dots \otimes e^{(N)}_{i_N}$ when $1 \leqslant i_k \leqslant n_k$.

Denote $\mathcal{B}=(f_1, \dots, f_p)$ the basis of $\mathcal{M}_1 \times \dots \times \mathcal{M}_N$ obtained from the bases $\mathcal{B}_i$. Then the family $\mathcal{F}'$ whose elements are the $\overline{f_{i_1} \otimes f_{i_2} \otimes \dots \otimes f_{i_N}}$ when $1 \leqslant i_1 \leqslant i_2 \leqslant \dots \leqslant i_N \leqslant n_1+\dots+n_N$ is a basis of $\Sym^N(\mathcal{M}_1 \otimes \dots \otimes \mathcal{M}_N)$ (cf \cite[p. 218]{Greub}). Moreover, $$\varphi(e^{(1)}_{i_1} \otimes e^{(2)}_{i_2} \otimes \dots \otimes e^{(N)}_{i_N})=\overline{f_{i_1} \otimes f_{n_1+i_2} \otimes \dots \otimes f_{n_1+\dots+n_{N-1}+i_N}}$$ so $\mathcal{F}$ is sent by $\varphi$ on a free family of $\Sym^N(\mathcal{M}_1 \otimes \dots \otimes \mathcal{M}_N)$, which proves that $\varphi$ is injective.

Hence the combination of \textbf{a)}, \textbf{b)} and \textbf{c)} in Proposition \ref{prop:tailleetopmodules} yields $$\sigma(\mathcal{M}_1 \otimes \dots \otimes \mathcal{M}_N) \leqslant \left(1+\log(N)\right)\max\big(\sigma(\mathcal{M}_1), \dots, \sigma(\mathcal{M}_N)\big).$$\end{dem}

\section[Chudnovsky's Theorem for Nilsson-Gevrey series]{A generalization of Chudnovsky's Theorem for Nilsson-Ge\-vrey series of arithmetic and holonomic type} \label{sec:chudnovskynilssongevrey}

The goal of this section is to prove Theorem \ref{th:chudnovskynilssongevreyandre}, a quantitative version of André's result on a Chudnovsky type Theorem for Nilsson-Gevrey series of arithmetic type \cite{AndregevreyI}.

In what follows, we will denote the standard derivation $\mathrm{d}/\mathrm{d}z$ on $\Qbar(z) $ by $\partial$.

\subsection{Product and sum of solutions of $G$-operators} 

The following result is proved by André in \cite[p. 720]{AndregevreyI}:

\begin{prop}[\cite{AndregevreyI}, p. 720] \label{prop:andreproduitsommegfonctions}
If $y_1$ (resp. $y_2$) is solution of a $G$-operator $L_1$ (resp. $L_2$), then
\begin{enumerateth}
\item $y_1+y_2$ is solution of a $G$-operator $L_3$;
\item $y_1 y_2$ is solution of a $G$-operator $L_4$.
\end{enumerateth}
\end{prop}

Note that $y_1$ and $y_2$ need not be $G$-functions, but they are Nilsson-Gevrey series of arithmetic and holonomic type. This proposition is obvious if $y_1$ and $y_2$ are $G$-functions, because the set of $G$-functions is a ring.

We are going to use the results of Section \ref{sec:sizediffmodule} (Proposition \ref{prop:tailleetopmodules}) in order to bound the sizes of $L_3$ and $L_4$.
\begin{prop}\label{prop:effectifproduitsommegfonctions}
Let $L_1$ and $L_2$ be the respective minimal operators of $y_1$ and $y_2$, which are $G$-operators. We can find nonzero $G$-operators $L_3$ and $L_4$ as in Proposition \ref{prop:andreproduitsommegfonctions}  such that $$\sigma(L_3) \leqslant \max\big(\sigma(L_1),\sigma(L_2)\big) \quad \text{and} \quad \sigma(L_4) \leqslant \big(1+\log(2)\big) \max\big(\sigma(L_1),\sigma(L_2)\big).$$
\end{prop}

\begin{rqu}
The result of Proposition \ref{prop:effectifproduitsommegfonctions} still holds for non minimal operators. Indeed, if $\widetilde{L}_i(y_i(z))=0$ then $\widetilde{L}_i$ is a multiple in $\Qbar(z)[\partial]$ of $L_i$ and we then have $\sigma(L_i) \leqslant \sigma(\widetilde{L}_i)$ (see Theorem  \ref{th:tailleproduitgop}).
\end{rqu}

\begin{dem}[of Propositions \ref{prop:andreproduitsommegfonctions} and \ref{prop:effectifproduitsommegfonctions}]
In \cite[p. 720]{Andre}, André gave a sketch of proof of Proposition~\ref{prop:andreproduitsommegfonctions}, which we make explicit here in order to bound the size of the operators.

If $L \in \Qbar(z)[\partial]$ is the minimal operator of a function $y$, we notice that there is a natural isomorphism
of differential modules between $\Qbar(z)[\partial]/\Qbar(z)[\partial]L$ and $$\Qbar(z)[\partial](y) :=\big\{M(y) \mid M \in \Qbar(z)[\partial]\big\}$$ given by $M \mod L \mapsto M(y)$. Hence $\sigma(L)=\sigma\left(\left(\Qbar(z)[\partial](y)\right)^*\right)$.

We also define for $u,v$ solutions of $G$-operators, $\Qbar(z)[\partial](u, v):= \left(\Qbar(z)[\partial](u)\right)[\partial](v)$ which is the set of linear combinations with coefficients in $\Qbar(z)$ of the $\partial^k(u) \partial^{\ell}(v)$, for $k, \ell \in \N$. The existence of differential equations over $\Qbar(z)$ satisfied by $u$ and $v$ ensures that this is indeed a finite dimensional $\Qbar(z)$-vector space, so that it can be endowed with a structure of differential module.

\begin{itemize}
\item Let $L_3$ be the minimal operator of $y_1+y_2$ over $\Qbar(z)$. Then $\sigma(L_3)=\sigma\left(\left(\Qbar(z)[\partial](y_1+y_2)\right)^*\right)$.

Let $\mathcal{M}$ be the image of the morphism of left $\Qbar[z][\partial]$-modules $$\fonction{\iota}{\Qbar(z)[\partial]}{\Qbar(z)[\partial](y_1) \times \Qbar(z)[\partial](y_2)}{L}{(L(y_1), L(y_2)).}$$ Thus, $\mathcal{M}$ is a differential submodule of $\Qbar(z)[\partial](y_1) \times \Qbar(z)[\partial](y_2)$ because it is a sub-left-$\Qbar(z)[\partial]$-module of $\Qbar(z)[\partial](y_1) \times \Qbar(z)[\partial](y_2)$ which is of finite dimension over $\Qbar(z)$.

Moreover, set $$\fonction{\varphi}{\mathcal{M}}{\Qbar(z)[\partial](y_1+y_2)}{(L(y_1), L(y_2))}{L(y_1+y_2)=L(y_1)+L(y_2).}$$ This morphism of differential modules is well-defined and surjective, hence its factorisation by its kernel induces an isomorphism of differential modules between $\Qbar(z)[\partial](y_1+y_2)$ and $\mathcal{M}/\ker(\varphi)$. Therefore, there is an isomorphism $\left(\Qbar(z)[\partial](y_1+y_2)\right)^* \simeq \left(\mathcal{M}/\ker(\varphi)\right)^* \simeq \mathcal{M}^*/\Img(\varphi^*)$ where $\varphi^*$ is the dual morphism of $\varphi$ in the sense of linear algebra. Hence by Proposition \ref{prop:tailleetopmodules} \textbf{a)}, we get \begin{equation}\label{eq:effectifandre1}
    \sigma(L_3)=\sigma(\Qbar(z)[\partial](y_1+y_2)^*) \leqslant \sigma(\mathcal{M}^*) \end{equation}

We now proceed to bound $\sigma(\mathcal{M}^*)$. The dual morphism of the injection $\mathcal{M} \hookrightarrow \Qbar(z)[\partial](y_1) \times \Qbar(z)[\partial](y_2)$ is a surjection $\left(\Qbar(z)[\partial](y_1)\times \Qbar(z)[\partial](y_2)\right)^* \twoheadrightarrow \mathcal{M}^*$, and $$\left(\Qbar(z)[\partial](y_1) \times \Qbar(z)[\partial](y_2)\right)^* \simeq \left(\Qbar(z)[\partial](y_1)\right)^* \times \left(\Qbar(z)[\partial](y_2)\right)^*$$ so that $\mathcal{M}^*$ can be written as some quotient of $\left(\Qbar(z)[\partial](y_1)\right)^* \times \left(\Qbar(z)[\partial](y_2)\right)^*$.

Thus, it follows from Proposition \ref{prop:tailleetopmodules} \textbf{a)} and \textbf{b)} that \begin{equation}\label{eq:effectifandre2} \sigma(\mathcal{M}^*) \leqslant \sigma\left( \left(\Qbar(z)[\partial](y_1)\right)^* \times \left(\Qbar(z)[\partial](y_2)\right)^*\right) \leqslant \max\big(\sigma(L_1),\sigma(L_2)\big).
\end{equation}

Finally, \eqref{eq:effectifandre1} and \eqref{eq:effectifandre2} imply that $\sigma(L_3) \leqslant \max\big(\sigma(L_1),\sigma(L_2)\big) < \infty.$ Hence $L_3$ is indeed a $G$-operator.

\item We define the morphism $$\fonction{\psi}{\Qbar(z)[\partial](y_1) \otimes \Qbar(z)[\partial](y_2)}{\Qbar(z)[\partial](y_1, y_2)}{K(y_1) \otimes L(y_2)}{K(y_1)L(y_2).}$$ Moreover, the set $\mathcal{M}=\Qbar(z)[\partial](y_1 y_2)$ is a differential submodule of $\Qbar(z)[\partial](y_1, y_2)$ because if $L \in \Qbar(z)[\partial]$, $L=\sum\limits_{k=0}^{\mu} a_k \partial^k$, we have $$L(y_1 y_2)=\sum_{k=0}^{\mu} a_k \sum_{p=0}^k \binom{k}{p} \partial^p(y_1) \partial^{k-p}(y_2)=\psi\left(\sum_{k=0}^{\mu} a_k \sum_{p=0}^k \binom{k}{p} \partial^p(y_1) \otimes \partial^{k-p}(y_2)\right) \in \Qbar(z)[\partial](y_1, y_2).$$ By definition, the restriction of $\psi$ to $\psi^{-1}(\mathcal{M})$ is a surjection $\widetilde{\psi} : \psi^{-1}(\mathcal{M}) \rightarrow \mathcal{M}$.

The factorisation of $\widetilde{\psi}$ by its kernel shows that $\mathcal{M}$ is isomorphic to $\mathcal{N}/\ker(\widetilde{\psi})$ where $\mathcal{N}=\psi^{-1}(\mathcal{M})$ is a submodule of $\Qbar(z)[\partial](y_1) \otimes \Qbar(z)[\partial](y_2)$. 

 Passing to the dual modules, we get $\mathcal{M}^* \simeq \mathcal{N}^*/\ker(\widetilde{\psi}^*)$ hence by Proposition \ref{prop:tailleetopmodules} \textbf{a)}, if $L_4$ is the minimal operator of $y_1 y_2$ over $\Qbar(z)$, then \begin{equation} \label{eq:effectifandre3}\sigma(L_4) =\sigma(\mathcal{M}^*) \leqslant \sigma(\mathcal{N}^*).\end{equation}

Furthermore, the dual morphism of $\mathcal{N} \hookrightarrow \Qbar(z)[\partial](y_1) \otimes \Qbar(z)[\partial](y_2)$ is a surjection $$\left(\Qbar(z)[\partial](y_1) \otimes \Qbar(z)[\partial](y_2)\right)^* \twoheadrightarrow \mathcal{N}^*$$ and we have $$\left(\Qbar(z)[\partial](y_1) \otimes \Qbar(z)[\partial](y_2)\right)^* \simeq \left(\Qbar(z)[\partial](y_1)\right)^* \otimes \left(\Qbar(z)[\partial](y_2)\right)^*\;,$$ so that, by Proposition \ref{prop:tailleetopmodules} \textbf{a)} and Proposition \ref{prop:tailleetproduittensoriel},
\begin{equation} \label{eq:effectifandre4}\sigma(\mathcal{N}^*) \leqslant \sigma\left(\left(\Qbar(z)[\partial](y_1)\right)^* \otimes \left(\Qbar(z)[\partial](y_2)\right)^*\right) \leqslant \big(1+\log(2)\big) \max\big(\sigma(L_1),\sigma(L_2)\big).\end{equation}
\end{itemize}

Thus, the combination of \eqref{eq:effectifandre3} and \eqref{eq:effectifandre4} shows that $\sigma(L_4) \leqslant \big(1+\log(2)\big) \max\big(\sigma(L_1),\sigma(L_2)\big)$, so that $L_4$ is a $G$-operator.\end{dem}

\begin{rqu}
Using Proposition \ref{prop:tailleetproduittensoriel}, we can generalize the second statement of Proposition \ref{prop:effectifproduitsommegfonctions} to an arbitrary product $y_1 \dots y_N$ of solutions of $G$-operators: if $L_0$ is the minimal operator of $y_1 \dots y_N$ and $L_i$ is the minimal operator of $y_i$, we obtain $$\sigma(L_0) \leqslant \big(1+\log(N)\big) \max\big(\sigma(L_1),\dots,\sigma(L_N)\big).$$
\end{rqu}

\subsection{Proof of Theorem \ref{th:chudnovskynilssongevreyandre}}

We can now deduce from Proposition \ref{prop:andreproduitsommegfonctions} and its quantitative version the proof of Theorem~\ref{th:chudnovskynilssongevreyandre} stated in the introduction.

\begin{dem}[of Theorem \ref{th:chudnovskynilssongevreyandre}]
\begin{itemize}
For all $(\alpha, k, \ell) \in S$, let $N_{\alpha,k,\ell} \neq 0$ be the minimal operator over $\Qbar(z)$ of $g_{\alpha,k,\ell}(z):=c_{\alpha, k,\ell} z^{\alpha} \log(z)^k f_{\alpha, k,\ell}(z)$, with $f_{\alpha,k,\ell}(z)$ a $G$-function, which exists since this function is a product of solutions of differential operators with coefficients in $\Qbar(z)$.
    \item Since $f(z)$ is the sum over $(\alpha,k,\ell) \in S$ of the functions $g_{\alpha,k,\ell}(z)$, Proposition \ref{prop:effectifproduitsommegfonctions} implies that the size of the minimal operator $L$ of $f(z)$ over $\Qbar(z)$ satisfies $$\sigma(L) \leqslant \max\limits_{(\alpha, k, \ell) \in S} \sigma(N_{\alpha, k,\ell}).$$
    \item Let $(\alpha, k, \ell) \in S$. The function $z^{\alpha} \log(z)^k c_{\alpha, k,\ell} f_{\alpha,k,\ell}(z)$ is a product of $k+2$ solutions of $G$-operators. The factors of this product are :
    \begin{itemizeth}[label=--]
     \item $z^{\alpha}$ of minimal operator $K_{\alpha}=\mathrm{d}/\mathrm{d}z-\alpha/z$ over $\Qbar(z)$;
     \item $\log(z)$ of minimal operator $T=z(\mathrm{d}/\mathrm{d}z)^2+\mathrm{d}/\mathrm{d}z$ over $\Qbar(z)$ ($k$ times);
     \item $c_{\alpha,k,\ell} f_{\alpha,k,\ell}(z)$ of minimal operator $L_{\alpha,k,\ell}$ over $\Qbar(z)$.
     \end{itemizeth}

 Thus, it follows from the remark after Proposition \ref{prop:effectifproduitsommegfonctions} that $$\sigma(N_{\alpha, k,\ell}) \leqslant \big(1+\log(k+2)\big) \max\left(\sigma(K_{\alpha}), \sigma(T), \sigma(L_{\alpha, k,\ell})\right)$$

    But we can compute directly the iterated matrix corresponding to the differential system $y'=\gf{\alpha}{z}y$: denoting $B_{\alpha}=\gf{\alpha}{z}$, we have, with the notations of Section \ref{subsec:Sizeandops}, $$\forall s \in \N, \quad \left(B_{\alpha}\right)_s = \gf{\alpha(\alpha-1) \dots (\alpha-s+1)}{z^s}$$ so that $\den(\alpha)^{2s} z^s \left(B_{\alpha}\right)_s \in \Z[z]$ (see \cite[Lemma 10 p. 334]{Lepetit2}) whence $\sigma(K_{\alpha}) \leqslant 2 \log(\den(\alpha))$. Moreover, $\sigma(T) \leqslant 1$, so we finally get \begin{align*}
    &\sigma(L) \leqslant \max\limits_{(\alpha, k, \ell) \in S}\left(\left(1+\log(k+2)\right)\max \left(2 \log(\den(\alpha)), 1, \sigma(L_{\alpha, k,\ell})\right)\right) \\
    &\leqslant \max\left(1+\log(\kappa+2),\; 2(1+\log(\kappa+2))\log\big(\max_{\alpha \in A} \den(\alpha)\big),  \max\limits_{(\alpha, k, \ell) \in S}\left((1+\log(k+2))\sigma(L_{\alpha, k,\ell})\right) \right).
    \end{align*}
\end{itemize}
 \noindent where $\kappa$ is the maximum of  the integers $k$ such that $(\alpha, k, \ell) \in S$ for some $(\alpha, \ell) \in \Q \times \N$ and $A=\{ \alpha \in \Q : \exists (k, \ell) \in \N^2, (\alpha, k, \ell) \in S \}$. This completes the proof of Theorem \ref{th:chudnovskynilssongevreyandre}.
\end{dem}

\begin{rqu}
The quantitative version of Chudnovsky's Theorem in \cite[p. 299]{Dwork} implies that
\begin{equation}\label{eq:chudnovskydworkI}\sigma(L_{\alpha,k,\ell}) \leqslant \left(5 \mu_{\alpha,k,\ell}^2 (\delta_{\alpha,k,\ell}+1)-1-(\mu_{\alpha,k,\ell} - 1)(\delta_{\alpha,k,\ell}+1)\right)\overline{\sigma}(f_{\alpha,k,\ell})\end{equation}
provided that $L_{\alpha,k,\ell} \in \K\left[z,\mathrm{d}/\mathrm{d}z\right]$, where $\K$ is a number field, and $\mu_{\alpha,k,\ell} :=\ord(L_{\alpha,k,\ell}) \geqslant 2$. In \eqref{eq:chudnovskydworkI}, we denote $\delta_{\alpha,k,\ell}:=\deg_z(L_{\alpha,k,\ell})$ and $\overline{\sigma}(y) :=\sigma(y)+\limsup\limits_{s \rightarrow + \infty} \gf{1}{s} \sum\limits_{\tau : \K \hookrightarrow \C} \sup\limits_{m \leqslant s} \log^{+}|y_m|_{\tau}$ where $\sigma(y)$ is defined in Definition \ref{def:taillematrice} above, and for any embedding $\tau : \K \hookrightarrow \C$, $$ |\zeta|_{\tau}:=\begin{cases} |\tau(\zeta)|^{1/[\K:\Q]} \quad \mathrm{if} \; \tau(\K) \subset \R \\ |\tau(\zeta)|^{2/[\K:\Q]} \quad \mathrm{else}. \end{cases}$$

On the other hand, if $\mu_{\alpha,k,\ell}=1$, we have \begin{equation}\label{eq:chudnovskydworkII}
    \sigma(L_{\alpha,k,\ell}) \leqslant (6\delta_{\alpha,k,\ell} - 1)\overline{\sigma}(f_{\alpha,k,\ell}).
\end{equation}

Therefore, the combination of Theorem \ref{th:chudnovskynilssongevreyandre} and of \eqref{eq:chudnovskydworkI} and \eqref{eq:chudnovskydworkII} (\emph{i.e} Chudnovsky's Theorem) gives an upper bound for $\sigma(L)$ in terms of the $\overline{\sigma}(f_{\alpha,k,\ell})$, for $(\alpha,k,\ell) \in S$.
\end{rqu}

\section{Applications} \label{sec:application}

We now give two applications of the results of Sections \ref{sec:sizediffmodule} and \ref{sec:chudnovskynilssongevrey}. The first one consists in expressing the size of a product of $G$-operators in terms of the sizes of each term of the product; the second one is related to a Diophantine problem studied in \cite{Lepetit2}.

\subsection{Size of a product of $G$-operators} \label{subsec:sizeproductGop}

In this part we are going to interpret André's result on the size of differential modules (Proposition \ref{prop:tailleetopmodules}) in terms of differential operators. 

The following proposition enables us to formulate a correspondance between the right factors of a differential operator $L$ and the differential submodules of $\mathcal{M}_L$ (see \cite[p. 47, p. 58]{Singer}). We denote by $\partial$ the derivation $\mathrm{d}/\mathrm{d}z$ on $\Qbar(z)$.

\begin{prop}[\cite{Singer}, p. 47] \label{prop:L1L2exactseq}
Given $L_1, L_2 \in \Qbar(z)[\partial]$, we have an exact sequence \begin{equation}\label{eq:seqexacteL1L2} 0 \rightarrow \mathcal{M}_{L_2} \rightarrow \mathcal{M}_{L_1 L_2} \rightarrow \mathcal{M}_{L_1} \rightarrow 0.\end{equation}
\end{prop}

\begin{dem}
Since we are working with finite-dimensional vector spaces, a sequence $$0 \rightarrow \mathcal{M}_{L_2} \xrightarrow{u} \mathcal{M}_{L_1 L_2} \xrightarrow{v} \mathcal{M}_{L_1} \rightarrow 0$$ is exact if and only if the dual one $$0 \leftarrow \mathcal{M}^*_{L_2} \xleftarrow{u^*} \mathcal{M}^*_{L_1 L_2} \xleftarrow{v^*} \mathcal{M}^*_{L_1} \leftarrow 0$$ is exact (see \cite[pp. 53 -- 59]{Jacobson}). Thus, it suffices to find an exact sequence \begin{equation}\label{eq:seqexacteL1L2II} 0 \rightarrow \Qbar(z)[\partial]/\Qbar(z)[\partial]L_1 \xrightarrow{\varphi} \Qbar(z)[\partial]/\Qbar(z)[\partial]L_1 L_2 \xrightarrow{\psi} \Qbar(z)[\partial]/\Qbar(z)[\partial]L_2 \rightarrow 0.\end{equation}

For this purpose, we define the injective map $$\fonction{\varphi}{\Qbar(z)[\partial]/\Qbar(z)[\partial]L_1}{\Qbar(z)[\partial]/\Qbar(z)[\partial]L_1L_2}{u \mod L_1}{uL_2 \mod L_1L_2}$$ and the surjective map $$\fonction{\psi}{\Qbar(z)[\partial]/\Qbar(z)[\partial]L_1 L_2}{\Qbar(z)[\partial]/\Qbar(z)[\partial]L_2}{u \mod L_1L_2}{u \mod L_2.}$$ We have $\ker\left(\psi\right)=\left.\Qbar(z)[\partial]L_2\middle/\Qbar(z)[\partial]L_1 L_2\right.=\Img\left(\varphi\right)$. Hence the sequence \eqref{eq:seqexacteL1L2II} is indeed an exact one.
\end{dem}

\begin{rqu}
A practical consequence of Proposition \ref{prop:L1L2exactseq} is the equivalence over $\Qbar(z)$ of the differential system $Y'=A_{L_1 L_2} Y$ to some system $Z'=BZ$, where $B=\begin{pmatrix} A_{L_1} & B^{0} \\ 0 & A_{L_2} \end{pmatrix}$ and $B^0$ has coefficients in $\Qbar(z)$. This follows from Lemma \ref{lem:formematriceseqexacte} applied to \eqref{eq:seqexacteL1L2}.
\end{rqu}

\bigskip

Propositions \ref{prop:tailleetopmodules} and \ref{prop:L1L2exactseq} then imply the following result, which is the desired application.

\begin{Th} \label{th:tailleproduitgop}
For any $L_1, L_2 \in \Qbar(z)[\partial]$, the following inequalities hold: \begin{align*} \max\big(\sigma(L_1), \sigma(L_2)\big) \leqslant \sigma(L_1 L_2) &\leqslant 1 + \gf{11}{6} \max\left(\sigma(L_1), \sigma(L_2), \sigma(L_1^*)\right) \\ &\leqslant 1 + \gf{11}{6} \max\big(\sigma(L_1)+\ord(L_1)-1, \sigma(L_2)\big)\end{align*} and $$\sigma(L_1 L_2) \leqslant \ord(L_1)+ 2 \sigma(L_1)+ \sigma(L_2).$$
\end{Th}

\begin{dem}
Proposition \ref{prop:L1L2exactseq} ensures us that the sequence $$ 0 \rightarrow \mathcal{M}_{L_2} \rightarrow \mathcal{M}_{L_1 L_2} \rightarrow \mathcal{M}_{L_1} \rightarrow 0$$ is exact. Therefore, $\mathcal{M}_{L_2}$ is a differential submodule of $\mathcal{M}_{L_1 L_2}$ and $\mathcal{M}_{L_1} \simeq \left.\mathcal{M}_{L_1 L_2} \middle/ \mathcal{M}_{L_2}\right.$.  It then follows from Proposition \ref{prop:tailleetopmodules} \textbf{a)} that $$\max(\sigma(L_1),\sigma(L_2))=\max\left(\sigma(\mathcal{M}_{L_1}),\sigma(\mathcal{M}_{L_2})\right) \leqslant \sigma(\mathcal{M}_{L_1L_2})=\sigma(L_1 L_2)$$ which is the first inequality we wanted to prove.

On the other hand, Equation \eqref{eq:seqexacteeq2} and Proposition \ref{prop:tailleetopmodules} \textbf{b)} yield $$\sigma(L_1 L_2) \leqslant 1+\gf{11}{6} \max\left(\sigma(L_1),\sigma(L_2),\sigma(L_1^*)\right).$$ Moreover, \eqref{eq:sizeadjalt} after Proposition \ref{prop:tailleetopmodules} implies $\sigma(L_1^*) \leqslant \sigma(L_1)+\ord L_1-1$, so that $$\sigma(L_1 L_2) \leqslant 1+\gf{11}{6}\max\big(\sigma(L_1)+\ord L_1-1, \sigma(L_2)\big).$$
Likewise, we obtain the inequality $\sigma(L_1 L_2) \leqslant \ord(L_1)+ 2 \sigma(L_1)+ \sigma(L_2)$ by applying Assertion~\eqref{eq:seqexacteeq3}.
\end{dem}

\subsection{Application to a Diophantine problem} \label{subsec:diophproblem}

Let us now consider a Diophantine approximation problem. It is a generalization of results by Fischler and Rivoal \cite{FRivoal} studied in \cite{Lepetit2}.

Let $\K$ be a number field and $F(z)=\sum\limits_{k=0}^{\infty} A_k z^k \in \K\llbracket z\rrbracket$ a nonpolynomial $G$-function of radius of convergence $R >0$. Let $L \in \K\left[z, \mathrm{d}/\mathrm{d}z \right] \setminus \{0\}$ the minimal operator of order $\mu$ of $F$, which is therefore a $G$-operator.  
Take a parameter $\beta \in \Q \setminus \Z_{\leqslant 0}$, that will remain fixed. For $n \in \N^*$ and $s \in \N$, we define the $G$-functions $$F_{\beta,n}^{[s]}(z)=\sum\limits_{k=0}^{\infty} \gf{A_k}{(k+\beta+n)^s} z^{k+n}.$$ These are related to iterated primitives of $F(z)$. The Diophantine problem we are interested in is to find upper and lower bounds on the dimension of the vector space $$\Phi_{\alpha, \beta, S} :=\Vect_{\K}\big(F_{\beta,n}^{[s]}(\alpha), \; n \in \N, \;\; 0 \leqslant s \leqslant S\big)$$ when $S$ is a large enough integer and $\alpha \in \K$, $0<|\alpha| <R$. Note that it is not obvious that $\Phi_{\alpha, \beta, S}$ has finite dimension. Precisely, in \cite{Lepetit2}, we prove the following theorem:

\begin{Th} \label{th:chapitre2}
Assume that $F$ is not a polynomial. Then for $S$ large enough, the following inequality holds: $$ \gf{1+o(1)}{[\K:\Q] C(F,\beta)} \log(S) \leqslant \dim_{\K} \Phi_{\alpha,\beta,S} \leqslant \ell_0(\beta) S+\mu.$$ Here, if $\delta=\deg_z(L)$ and $\omega$ is the order of $0$ as a singularity of $L$, $\ell_0(\beta)$ is defined as the maximum of $\ell:=\delta-\omega$ and the numbers $f-\beta$ when $f$ runs through the exponents of $L$ at infinity such that $f-\beta \in \N$, and $C(F,\beta)$ is a positive constant depending only on $F$ and $\beta$, and not $\alpha$.
\end{Th}

Fischler and Rivoal proved this theorem for $\beta=0$ in \cite{FRivoal}. With their method, the constant $C(F,0)$ is computable in principle, but they didn't give an explicit formula for it. In \cite{Lepetit2}, we express explicitly $C(F,\beta)$ in function of quantities depending on $L$, $F$ and the denominator of $\beta$ using the theory of $G$-operators. Let us explain how Theorems \ref{th:chudnovskynilssongevreyandreI} and \ref{th:tailleproduitgop} of the present paper can be used to compute $C(F,\beta)$.

By \cite[Lemma 1 p. 11]{FRivoal}, we can find polynomials $Q_0(X), \dots, Q_{\ell}(X) \in \Oal_{\K}[X]$ and $u\in \N^*$ such that $$u z^{\mu-\omega} L=\sum_{j=0}^{\ell} z^j Q_j(\theta+j),$$ with $\theta=z\mathrm{d}/\mathrm{d}z$, $\mu$ the order of $L$, $\omega$ the multiplicity of $0$ as a singularity of $L$ and $\ell=\delta-\omega$ where $\delta$ is the degree in $z$ of $L$. We can show that if $\ell=0$, then $F(z) \in \K[z]$, so that this case is excluded by the assumption of Theorem \ref{th:chapitre2}.

Define, for $j \in \{ 0, \dots, \ell \}$, $Q_{j, \beta}(X) := Q_j(X-\beta)$. Then the differential operator of order $\mu$ \begin{equation} \label{eq:defLbeta}L_{\beta} :=\sum\limits_{j=0}^{\ell} z^j Q_{j,\beta}(\theta+j) \in \K\left[z,\gf{\mathrm{d}}{\mathrm{d}z}\right]\end{equation} is the minimal operator over $\Qbar(z)$ of the Nilsson-Gevrey series of arithmetic and holonomic type $z^{\beta} F(z)$. We introduce the operator $$\widetilde{L}_{\beta}=\left(\gf{\mathrm{d}}{\mathrm{d}z}\right)^{\ell} z^{m-1} L_{\beta}.$$ In \cite{Lepetit2}, a crucial point, which we couldn't solve without using the results of the present paper, is the evaluation of  $\sigma(\widetilde{L}_{\beta})$ in function of $\sigma(L)$ or of $\overline{\sigma}(F)$. Indeed, $\sigma(\widetilde{L}_{\beta})$ occur in the expression of $C(F,\beta)$.

We can now answer that question. Indeed, on the one hand, if $L_0=\left(\mathrm{d}/\mathrm{d}z\right)^{\ell}$, Theorem \ref{th:tailleproduitgop} implies that $$\sigma(\widetilde{L}_{\beta}) \leqslant \ell + 2 \sigma(L_0) + \sigma(L_{\beta})$$ since $\sigma(z^{m-1} L_{\beta})=\sigma(L_{\beta})$. Moreover, a basis of solutions of the equation $L_0(y(z))=0$ is $(1, z, \dots, z^{\ell-1})$ so that a fundamental matrix of solution of the system $y'=A_{L_0} y$ is the wronskian matrix $$Y=\begin{pmatrix} 1 & z & \dots & z^{\ell-1} \\ 0 & 1 & & (\ell-1) z^{\ell-2} \\ \vdots & 0 & \ddots & \vdots \\ 0 & 0 & \dots & (\ell-1)! \end{pmatrix},$$ which satisfies $Y^{(s)}=0$ for $s$ large enough. Hence $(A_{L_0})_s=Y_s Y^{-1}=0$ for $s$ large enough and $\sigma(L_0)=0$.

On the other hand, by applying Theorem \ref{th:chudnovskynilssongevreyandre} to the Nilsson-Gevrey series $z^{\beta} F(z)$, we obtain \begin{equation}\sigma(L_{\beta}) \leqslant \big(1+\log(2)\big) \max\left(1, 2 \log(\den(\beta)), \sigma(L) \right).\end{equation} Finally, using the bound on $\sigma(L)$ arising from Chudnovsky's Theorem (\eqref{eq:chudnovskydworkI} and \eqref{eq:chudnovskydworkII}), we have

\begin{equation}\label{eq:finaleapplicationchud}\sigma(L_{\beta}) \leqslant
 \big(1+\log(2)\big) \max\big(1, 2 \log(\den(\beta)),\left((5+\varepsilon_{\mu,1}) \mu^2 (\delta+1)-1-(\mu - 1)(\delta+1)\right)\overline{\sigma}(F) \big),  \end{equation} where $\varepsilon_{x,y}=\begin{cases} 1 & \text{if} \; x=y \\ 0 & \text{else} \end{cases}$ denotes the Kronecker symbol, and $$\sigma(\widetilde{L}_{\beta}) \leqslant \ell + \big(1+\log(2)\big) \max\big(1, 2 \log(\den(\beta)),\left((5+\varepsilon_{\mu,1}) \mu^2 (\delta+1)-1-(\mu - 1)(\delta+1)\right)\overline{\sigma}(F) \big).$$ This is the desired bound.

\begin{rqu}
In \cite{LepetitNilssonG}, we tried another method to obtain an inequality of the type of \eqref{eq:finaleapplicationchud}: we adapted directly the proof of Chudnovsky's Theorem to the case of a subclass of the set of Nilsson-Gevrey series of arithmetic and holonomic type. We now explain this approach.  

Let $\K$ be a number field and \begin{equation}\label{eq:NGsubclass} f(z)=\sum\limits_{\ell=1}^{\lambda} \sum\limits_{k=0}^{\kappa} z^{\alpha_{\ell}} \log(z)^k f_{k,\ell}(z) \end{equation} be a Nilsson-Gevrey series of arithmetic and holonomic type of order $0$, where $\boldsymbol{\alpha}=(\alpha_1, \dots, \alpha_{\lambda}) \in \Q^{\lambda}$, and the $f_{k,\ell}(z) \in \K\llbracket z \rrbracket$ are $G$-functions. Here, $f(z)$ is a particular case of Nilsson-Gevrey series of arithmetic and holonomic type, which is not of the most general form \eqref{eq:defNGseries} given in the introduction. Indeed, in \eqref{eq:defNGseries}, it is allowed to take $f_{k,\ell}$ as a linear combination with coefficients in $\C$ of $G$-functions, which is not permitted here. Then Proposition \ref{prop:chudNGalt} below provides an estimate on the size of the minimal operator of $f(z)$ over $\Qbar(z)$. This result uses a stronger hypothesis than Theorem \ref{th:chudnovskynilssongevreyandre}, namely that the minimal operator of $f$ over the fraction field $\mathcal{R}$ of the ring of "$\mathcal{S}$-polynomials" is the same as the minimal operator of $f$ over $\K(z)$.The $\mathcal{S}$-polynomials are defined as the functions of the form \eqref{eq:NGsubclass}, where $(\alpha_1, \dots, \alpha_{\lambda}) \in \Q^{\lambda}$ and the $f_{k,\ell}(z) \in \K[z]$.

\begin{prop} \label{prop:chudNGalt}
Let $f(z)$ be a function of the form \eqref{eq:NGsubclass}. \emph{Assume that the minimal operator $L$ of $f(z)$ of order $\mu$ over the fraction field $\mathcal{R}$ of the ring of the $\mathcal{S}$-polynomials satisfies $L \in \K(z)\left[\mathrm{d}/\mathrm{d}z\right]$}. Then we have  \begin{equation}\label{eq:chudNGAlteq1}\limsup_{s \rightarrow +\infty} \gf{1}{s} \log(q_s) \leqslant [\K:\Q]\left(2 \log(\den(\boldsymbol{\alpha}))+\kappa+ 2 \mu \left(2 \mu \lambda(\kappa+1)+1\right)(\delta+1) \max_{0 \leqslant k \leqslant \kappa \atop 1 \leqslant \ell \leqslant \lambda} \overline{\sigma}(f_{k,\ell})\right)\end{equation} with the notations of Definitions \ref{def:galochkin} and \ref{def:taillematrice} (with $G$ the companion matrix of $L$), where $\delta=\deg_z(L)$. 
\end{prop}

Let us compare the conclusion of Proposition \ref{prop:chudNGalt} with the one of Theorem \ref{th:chudnovskynilssongevreyandre} in the particular case of $L_{\beta}$. We can prove that $L_{\beta}$ is at the same time the minimal operator of $z^{\beta} F(z)$ over $\Qbar(z)$ and over $\mathcal{R}$.

If we consider the family $(f_1, \dots, f_{\mu}):=(F,F', \dots, F^{(\mu-1)})$ consisting of linearly independent functions over $\Qbar(z)$  and $(g_1, \dots, g_{\mu}):=(G,G', \dots, G^{(\mu-1)})$, where $G(z)=z^{\beta}F(z)$, then we see that $L_{\beta}=z^{\beta} L z^{-\beta}$ is the minimal operator of $G(z)$ over $\Qbar(z)$, so that $(g_1, \dots, g_{\mu})$ is free over $\Qbar(z)$.

The Leibniz formula then shows that $(g_1, \dots, g_{\mu}) \in z^{\beta-\mu} \K\llbracket z \rrbracket$. This implies that the family $(z^{-\beta+\mu} g_1, \dots, z^{-\beta+\mu} g_{\mu})$ is free over $\mathcal{R}$, using the linear independence over $\K\llbracket z \rrbracket$ of the family $(z^{\gamma} \log(z)^k)_{\gamma \in \Q \cap [0,1[ \atop k \in \N}$. Consequently, $(g_1, \dots, g_{\mu})$ is free over $\mathcal{R}$.

Thus, since the operator $L_{\beta}$ of order $\mu$ with coefficients in $\mathcal{R}$ annihilates $G(z)$, it is the minimal operator of $G(z)$ over $\mathcal{R}$.

In that case, using Proposition \ref{prop:lienq'ssigmaG}, Equation \eqref{eq:chudNGAlteq1} yields
   \begin{equation}\label{eq:conclusionChudNilssongevreyAlternatif}\sigma(L_{\beta}) \leqslant [\K:\Q]\left( 2\log(\den(\beta))+2 \mu(2\mu+1)(\delta+1) \overline{\sigma}(F)\right).\end{equation}
    \end{rqu}

To simplify, we assume that $\K=\Q$ and $d(\beta) \geqslant 2$ (\emph{i.e.} $\beta \not\in \Z$). 
\begin{itemize}[label=\textbullet]
    \item If $\den(\beta) \geqslant \exp\big((5+\varepsilon_{\mu,1}) \mu^2(\delta+1)-1-(\mu-1)(\delta+1))\overline{\sigma}(F)/2\big)$, then \eqref{eq:finaleapplicationchud} yields $$\sigma(L_{\beta}) \leqslant 2\big(1+\log(2)\big)\log(\den(\beta)), $$ which is a better estimate on $\sigma(L_{\beta})$ than \eqref{eq:conclusionChudNilssongevreyAlternatif}  if and only if $$\overline{\sigma}(F) > \gf{\log\big(\den(\beta)+2\big)}{\mu(2\mu+1)(\delta+1)}.$$
    
    \item Else, we see that the function $$\fonction{\mathrm{Comp}}{\N^* \times \N}{\R}{(\mu,\delta)}{ (1+\log(2))\big((5+\varepsilon_{\mu,1}) \mu^2 (\delta+1)-1-(\mu-1)(\delta+1)\big)-2\mu(2\mu+1)(\delta+1)}$$ is strictly positive. Indeed, if $\mu \geqslant 2$, $\mathrm{Comp}(\mu,\delta)=A(\mu)(\delta+1)-(1+\log(2))$ where $$A(\mu):=(3+5\log(2))\mu^2-(3+\log(2))\mu+1+\log(2)$$ is strictly greater than $1+\log(2)$ for all $\mu \geqslant 2$. Moreover, $\mathrm{Comp}(1,\delta)= 6 \log(2)(\delta+1)-1-\log(2) \geqslant 0$ for all $\delta \in \N$.
    
    Hence, if $(\mu,\delta) \in \N^* \times \N$ and $\den(\beta) \leqslant \exp\big(\overline{\sigma}(F) \mathrm{Comp}(\mu,\delta)/2\big)$, then \eqref{eq:conclusionChudNilssongevreyAlternatif} is a better estimate than \eqref{eq:finaleapplicationchud}.
    
    \item Else, if $\den(\beta)>\exp\big(\overline{\sigma}(F) \mathrm{Comp}(\mu,\delta)/2\big)$, then \eqref{eq:finaleapplicationchud} is a better estimate than \eqref{eq:conclusionChudNilssongevreyAlternatif}.
\end{itemize}

\bigskip

Let us finally consider the following explicit example: let $F(z)$ be the hypergeometric $G$-function $$F(z)= {}_2F_1\left(\gf{1}{3}, \gf{2}{11}, \gf{1}{6} ; z\right)=\sum_{k=0}^{\infty} \gf{\left(\frac{1}{3}\right)_k \left(\frac{2}{11}\right)_k}{\left(\frac{1}{6}\right)_k k!} z^k$$ whose minimal operator over $\Qbar(z)$ is $$L=z(z-1)\left(\ddz\right)^2+\left(\left(\gf{1}{3}+\gf{2}{11}+1\right)z-\gf{1}{6}\right) \left(\ddz\right)+\gf{1}{3} \times \gf{2}{11}.$$ We take $\beta=1/7$. It satisfies $\delta=\mu=2$. Hence we obtain with \eqref{eq:finaleapplicationchud} $\sigma(L_{\beta}) \leqslant 1232$ and with \eqref{eq:conclusionChudNilssongevreyAlternatif} $\sigma(L_{\beta}) \leqslant 784$. We see that our alternative method improves the bound \eqref{eq:finaleapplicationchud} on $\sigma(L_{\beta})$. 


\printbibliography

\bigskip
G. Lepetit, Université Grenoble Alpes, CNRS, Institut Fourier, 38000 Grenoble, France.

\url{gabriel.lepetit@univ-grenoble-alpes.fr}. 

\bigskip

\emph{Keywords}: $G$-functions, $G$-operators, Nilsson-Gevrey series of arithmetic type.

\bigskip

\emph{2020 Mathematics Subject Classification}. Primary 34M03, Secondary 13N10, 34M05, 11J72.
\end{document}